\documentclass[a4paper,twoside]{article}

\usepackage{graphs}
\usepackage{epsfig} % Mimi added this. Should we put it in graphs???

\author{Mireille Boutin and Gregor Kemper}

\title{Lossless Representation of Graphs using Distributions}

\renewcommand{\ss}[1]{_{\{#1\}}}
\newcommand{\ps}[1]{\phi\left(\{#1\}\right)}
\newcommand{\sps}[1]{_{\ps{#1}}}
\newcommand{\N}{\{1 \upto n\}}

\begin{document}

\maketitle

\begin{abstract}
We consider complete graphs with edge weights and/or node weights taking values in some set.
In the first part of this paper, we show that a large number of graphs 
are completely determined, up to isomorphism, by the distribution of their sub-triangles.
In the second part, we propose graph representations in terms of one-dimensional distributions
(e.g., distribution of the node weights, sum of adjacent weights, etc.).
For the case when the weights of the graph are real-valued vectors, 
we show that all graphs, except for a set of measure zero, are uniquely determined, 
up to isomorphism, from these distributions. 
The motivating application for this paper is the problem of browsing through large sets of graphs.
\end{abstract}

%\tableofcontents

\section*{Introduction} \label{0sIntro}

Graphs in general, 
and more particularly weighted or vector weighted graphs, 
are often used to represent complex structures such as 3D objects \cite{Foulds}.
This is often done, for example, when trying to determine whether two objects have a similar structure: 
by using a graph representation, 
the problem is simplified into determining whether the two underlying graphs are similar.
The problem of determining whether two graphs are the same, 
up to a relabeling of the nodes, 
is called the {\em graph isomorphism problem}.
While graph isomorphism is not an NP hard problem, it is still very hard.
In fact, the problem is sometimes assigned to a special complexity class called 
{\em graph isomorphism complete} (\mycite{Skiena}).

Graph isomorphism has been a mainstream research problem for more than 30 years (\mycite{Gati,ReadCorneil}).
Over these years, several solutions have been developed,
including graph matching methods 
(\mycite{CorneilGotlieb,Ullmann,KitchenRosenfeld,KimKim,FalkenhainerForbusGentner,HoraudSkordas,MyaengLopez}),
canonical labeling representation methods (\mycite{Nauty}),
graphs invariants (\mycite{CorneilKirkpatrick,Umeyama88,Chung,MessmerBunke99,MessmerBunke00}),
graph matching based on error-correcting isomorphism methods
(\mycite{TsaiFu,ShapiroHaralick,SanfeliuFu,EsheraFu,Wong}),
and approximate graph matching methods
(\mycite{HeraultHoraudVeillonNiez,KittlerChristmasPetrou,ChristmasKittlerPetrou,AlmohamadDuffuaa}).
In this paper, we show that many graphs can be uniquely
reconstructed from some simple distributions.
For example, in the first section of the paper,
we show that a large number of weighted graphs can be
uniquely reconstructed from their distribution of subtriangles.
In other words, their distribution of subtriangles provides a faithful (i.e., lossless) 
representation for these graphs.
This means that isomorphism can be detected
simply by comparison of the respective distribution of triangles.
For simpler comparison and visualization, we also introduce graph representation in terms of
one-dimensional distributions (e.g., node weights, edge weights and sum of adjacent edge weights).
As we show in Section \ref{3s:Mimi}, for many graphs, these representations are lossless. 
Actually, when the weights of the graphs take values inside $\RR$, 
the set of exception form a set of measure zero.
This work can be viewed as an extension to graphs
 of the representations we proposed in \cite{Boutin.Kemper} for point configurations.

The results we present are interesting both from a graph theory point point of view 
and from an application point of view. 
From an application point of view, 
having a faithful representation which can be compared quickly is useful in the case where
many graphs need to be compared in a short amount of time.
For example, an important problem is the problem of browsing for graphs within a large database.
In this problem, being able to compare graphs quickly is a key issue.
However, this is not the only issue. 
Indeed, exhaustive searches are not efficient ways to query a large database
because their complexity is proportional to the size of the database.
For faster search, an index structure must be built.
Database indices exploit the presence of natural clusters in the dataset
to successively rule out large regions of the data space.  
Unfortunately, 
some datasets do not exhibit natural clusters,
particularly high-dimensional ones (\mycite{BeyerGoldsteinRamakrishnanShaft99}).
Indeed, while some high-dimensional datasets can be dealt with successfully (\mycite{ShaftRamakrishnan05})
(e.g., when the underlying dimensionality of the dataset 
is much lower than the dimension of the vectors),
clusters can usually only be found in low-dimensional projections 
of the high-dimensional vectors (\mycite{WangYang05}).
However, projecting graph representations which are not invariant under isomorphism
is unproductive
because the labeling ambiguity blurs the distinction between the dimensions being projected.
Thus an isomorphism invariant representation is needed.
Moreover, a faithful representation guarantees the highest level of comparison accuracy. 
In particular, small distinguishing features are guaranteed to be preserved. 
In other words, for a large set of graphs,
the proposed representations  provide explicit coordinates to represent graphs,
and thus naturally lend themselves 
to the array of database projection and indexing techniques available in the literature.
This is in contrast with many current graph indexing approaches which first approximate 
the structure with a lossy representation before indexing.

From a graph theory point of view, 
there is a long tradition  of considering subgraphs of
a given graph and asking how much information about the graph is
contained in its subgraph structure. Perhaps the most important
example of this approach is Ulam's conjecture (see \mycite{Ulam}),
also called the reconstruction problem. The conjecture can be stated
as follows: Let $G$ be a simple non-directed graph with $n \ge 3$
nodes ({\em simple} means for each pair of nodes there either is an
edge between them or not). Take the set of all isomorphism classes of
$(n-1)$-subgraphs of $G$, i.e., all graphs obtained from $G$ by
deleting one node. To each of these isomorphism classes assign its
{\em multiplicity}, i.e., the number of nodes of $G$ whose deletion
leads to a graph in this isomorphism class. The set of isomorphism
classes together with their multiplicities is an example of a {\em
  multiset}. We call it the {\em distribution of $(n-1)$-subgraphs of
  $G$}. Ulam's conjecture states that every simple graph with $n \ge
3$ nodes is uniquely determined, up to isomorphism, by the
distribution of its $(n-1)$-subgraphs. This is often phrased by saying
that $G$ is reconstructible from its $(n-1)$-subgraphs. Ulam's
conjecture is still open. Using computational techniques,
\mycite{McKay97} showed that it holds for $n \le 11$.  See also
\mycite{Pouzet.Thiery} for some background and related questions.

Given the difficulty of Ulam's conjecture, there is no way that we
could expect that all non-isomorphic graphs can be distinguished by
their distributions of subtriangles. Indeed, the first counter example
occurs when one looks at graphs with~5 nodes. The example is given by
the two graphs in Figure~\ref{1fCounter}. 
So the question is just how well the distribution of subtriangles
discriminates graphs, and how it compares with other graph invariants.
As we show in Section \ref{1sReconstruct}, 
quite a large number of graphs can be uniquely represented from their distribution of triangles.

\providecommand{\point}{\circle*{2.5}}
\Figure{
  \fbox{
    \begin{picture}(90,90)
      \put(15,15){\point}
      \put(18,15){\line(1,0){24}}
      \put(45,15){\point}
      \put(48,15){\line(1,0){24}}
      \put(45,18){\line(0,1){24}}
      \put(75,15){\point}
      \put(45,45){\point}
      \put(45,48){\line(0,1){24}}
      \put(45,75){\point}
    \end{picture}
  }
  \hspace{11mm}
  \fbox{
    \begin{picture}(90,90)
      \put(15,15){\point}
      \put(18,15){\line(1,0){54}}
      \put(15,18){\line(0,1){54}}
      \put(75,15){\point}
      \put(75,18){\line(0,1){54}}
      \put(15,75){\point}
      \put(18,75){\line(1,0){54}}
      \put(75,75){\point}
      \put(45,45){\point}
    \end{picture}
  }
%  \hspace{8mm}
%  \fbox{
%    \begin{picture}(90,90)
%      \put(15,15){\point}
%      \put(75,15){\point}
%      \put(15,75){\point}
%      \put(75,75){\point}
%      \put(45,45){\point}
%      \put(48,48){\line(1,1){24}}
%      \put(42,42){\line(-1,-1){24}}
%      \put(42,48){\line(-1,1){24}}
%      \put(48,42){\line(1,-1){24}}
%    \end{picture}
%  }
} {Non-isomorphic graphs with the same distribution of subtriangles.}
{1fCounter}
%%% Local Variables: 
%%% mode: latex
%%% TeX-master: "graphs"
%%% End: 

\section{Preliminaries} \label{0sPreliminaries}

By a \df{graph} $G$ with~$n$ nodes we understand a function which
assigns to each subset $\{i,j\} \subseteq \N$ an element $g\ss{i,j}$
in some set $X$. The value $g\ss{i,j}$ is interpreted as the weight of
the edge between the nodes~$i$ and~$j$. Typically, $X$ will be $\RR$
or $\RR^d$, but any other set, including finite sets, are possible. In
other words, the graphs we consider are complete, undirected graphs
with weighted edges, including edges between a node and itself. The
values $g\ss{i,i}$ between a node and itself may be interpreted as a
node weight. Other classes of graphs are included as special cases:
\begin{itemize}
\item We may consider nodes~$i$ and~$j$ with $g\ss{i,j} = 0$ (or some
  other designated value in $X$) as not connected. So the case of
  incomplete graphs is included. The set of edges will be
  \begin{equation} \label{0eqEdges}
     \left\{\{i,j\} \subseteq \N \mid g\ss{i,j} \ne 0\right\}.
  \end{equation}
\item If all node weights $g\ss{i,j}$ are equal, this amounts to the
  same as saying that there are no node weights. So the case of graphs
  without any node weights and without edges between a node and itself
  is included. Likewise, we may consider graphs which have node
  weights but no edge weight.
\item If all $g\ss{i,j}$ lie in $\{0,1\}$, then $G$ may be considered
  as a simple graph, i.e., a graph which has no edge weights. The set
  of edges of $G$ are then given by~\eqref{0eqEdges}.
\end{itemize}

Let $G$ and $G'$ two graphs, both with~$n$ nodes and edge weights
$g\ss{i,j}$ and $g'\ss{i,j}$ in the same set $X$. We say that $G$ and
$G'$ are \df{isomorphic} if there exists a permutation $\pi \in S_n$
such that
\[
g'\ss{i,j} = g\ss{\pi(i),\pi(j)} \quad \text{for all} \quad i,j \in
\N.
\]
We write $G \cong G'$ for this.

Let $\mathcal D$ be a function which assigns a value to each graph $G$
with~$n$ nodes and edge weights in a set $X$. We say that $\mathcal D$
is a \df{class function} if $G \cong G'$ implies ${\mathcal D}(G) =
{\mathcal D}(G')$. Assume that $\mathcal D$ is a class function and
$G$ is a graph. We say that $G$ is \df{reconstructible from $\mathcal
  D$} (or, by abuse of language, reconstructible from ${\mathcal
  D}(G)$) if for all other graphs $G'$ with~$n$ nodes and edge weights
in the same set $X$, we have the implication
\[
{\mathcal D}(G') = {\mathcal D}(G) \quad \Longrightarrow \quad G'
\cong G.
\]
Note that saying that  $G$ is reconstructible from $\mathcal D$ is equivalent
to saying that $\mathcal D$ is a lossless representation of $G$.
Sometimes we will restrict $G$ and $G'$ to special classes of graphs,
such as graphs with edge weights but no node weights. So
reconstructibility of $G$ from $\mathcal D$ means that ${\mathcal
  D}(G)$ determines $G$ up to isomorphism. The word
``reconstructible'' should not be misunderstood to mean that we
actually have a method for constructing $G$, or an equivalent graph,
from the knowledge of ${\mathcal D}(G)$. In \sref{3s:Mimi} we will
consider various class functions ${\mathcal D}_1 \upto {\mathcal
  D}_r$. We say that $G$ is reconstructible from ${\mathcal D}_1 \upto
{\mathcal D}_r$ if for all other graphs $G'$ with~$n$ nodes and edge
weights in the same set $X$, we have the implication
\[
{\mathcal D}_i(G') = {\mathcal D}_i(G) \quad \text{for all} \quad i
\in \{1 \upto r\} \quad \Longrightarrow \quad G' \cong G.
\]

All functions $\mathcal D$ on graphs that we consider will assign a
distribution to a graph $G$, the simplest example being the
distribution of all edge weights $g\ss{i,j}$. To make this precise, we
use the concept of a multiset. A \df{multiset} is a set $\mathcal M$
together with a function ${\mathcal M} \to \NN_{>0}$ into the set of
positive integers, which we interpret as assigning a multiplicity to
each element of $\mathcal M$. Two multisets are considered equal if
their underlying sets and multiplicity functions are equal. If we
speak of the multiset consisting of all $a_i$ for some range of~$i$,
we mean that the multiplicity of each $a_i$ is $|\{j \mid a_j =
a_i\}|$.

A typical statement using this language would be that a graph $G$ with
edge weights $g\ss{i,j}$ is reconstructible from the distribution of
its edge weights $g\ss{i,j}$. It is easy to see that this statement is
true if and only if all edge weights of $G$ are equal. However, if we
restrict to graphs which have no node weights, then every graph with
at most three nodes is reconstructible from the distribution of its
edge weights.

%%% Local Variables: 
%%% mode: latex
%%% TeX-master: "graphs"
%%% End: 

\section{Reconstructibility from the distribution of subtriangles}
\label{1sReconstruct}

%However, in this paper we shift the focus by considering weighted
%graphs, so we assign a weight (e.g. a real number) to each edge.
%Simple graphs form a special case, because they occur as weighted
%graphs with weights~1 and~0, say. By making certain assumptions about
%the distinctness of the weights, we move away rather far from the case
%of simple graphs, but, as we will see shortly, this allows us to move
%into a region where reconstructibility from subtriangles can indeed be
%guaranteed.

Let $G$ be graph with $n \ge 3$ nodes and edge weights $g\ss{i,j} \in
X$. In order to get a more distinctive notation for edges between a
node and itself, we write $w_i := g\ss{i,i}$, which we interpret as
the weight of the node~$i$. For $i,j,k \in \N$ pairwise distinct, let
${\mathfrak t}\ss{i,j,k}$ be the multiset consisting of the ordered
pairs $(g\ss{i,j},w_k)$, $(g\ss{i,k},w_j)$, and $(g\ss{j,k},w_i) \in X
\times X$. So ${\mathfrak t}\ss{i,j,k}$ represents the subgraph of $G$
with nodes~$i$, $j$, $k$. Taking ${\mathfrak t}\ss{i,j,k}$ as a
multiset means that we do not assign any ordering on the nodes, so the
subgraph is considered up to isomorphism. Such a subgraph ${\mathfrak
  t}\ss{i,j,k}$ will be called a subtriangle. Moreover, let ${\mathcal
  T}_G$ be the multiset consisting of all ${\mathfrak t}_S$ for $S
\subseteq \N$ with $|S| = 3$. ${\mathcal T}_G$ is called the
\df{distribution of subtriangles} of $G$. If $G'$ is a graph with
weighted edges such that $G'$ is isomorphic to $G$, then clearly
${\mathcal T}_{G'} = {\mathcal T}_G$, so ${\mathcal T}$ is a class
function.

\subsection{A partial isomorphism} \label{1sPartial}

The following lemma establishes a ``partial isomorphism'' between
graphs having the same distribution of subtriangles.
\tref{1tReconstruct} and \cref{1cReconstruct} deal with situations
where this is a ``full'' isomorphism.

\begin{lemma} \label{1lReconstruct}
  In the above situation write
  \[
  P := \left\{\{i,j\} \mid 1 \le i < j \le n\right\}
  \]
  and
  \[
  E := \left\{S \in P \mid g_T \ne g_S \ \text{for all} \ T \in P
    \setminus \{S\}\right\}.
  \]
  Let $G'$ be a further graph with edge weights $g'\ss{i,j}$ and node
  weights $w'_i$ such that ${\mathcal T}_{G'} = {\mathcal T}_G$. Then
  there exists $\pi \in S_n$ (the set of bijective maps $\N \to \N$)
  such that
  \begin{enumerate}
  \item for all $i \in \N$ we have
    \[
    w_i = w'_{\pi(i)},
    \]
  \item for all $\{i,j\} \in E$ we have
    \[
    g\ss{i,j} = g'\ss{\pi(i),\pi(j)},
    \]
  \item for $i,j,k \in \N$ pairwise distinct with $\{i,k\} \in E$ and
    $\{j,k\} \in E$ we have
    \[
    g\ss{i,j} = g'\ss{\pi(i),\pi(j)}.
    \]
  \end{enumerate}
\end{lemma}

\begin{proof}
  Since ${\mathcal T}_G$ has $\binom{n}{3}$ elements, it follows that
  $G'$, like $G$, has~$n$ nodes.  For $S,T \in P$ we have by the
  definition of $E$:
  \begin{equation} \label{1eqInjective}
    \text{if} \quad S \in E \quad \text{or} \quad T \in E \quad
    \text{then} \quad g_S = g_T \quad \text{implies} \quad S = T.
  \end{equation}
  Form the union (with adding multiplicities) of all multisets
  ${\mathfrak t}\ss{i,j,k}$ lying in ${\mathcal T}_G$, and then take
  the multiset consisting of the first components of all
  $(g\ss{i,j},w_k)$ lying in this union. This yields the multiset
  consisting of all $g_S$, $S \in P$, counted $n - 2$ times for each
  $S \in P$. Since ${\mathcal T}_{G'} = {\mathcal T}_G$, this implies
  that the multiset of all $g_S$, $S \in P$ and the multiset of all
  $g'_S$, $S \in P$, coincide. With
  \[
  E' := \left\{S \in P \mid g'_T \ne g'_S \ \text{for all} \ T \in P
    \setminus \{S\}\right\},
  \]
  it follows that there exists a bijection $\map{\phi}{E}{E'}$ such
  that:
  \begin{equation} \label{1eqPhi}
    \text{for all} \quad S \in E \quad \text{we have} \quad
    g'_{\phi(S)} = g_S.
  \end{equation}
  From the definition of $E'$, we obtain for all $S,T \in P$:
  \begin{equation} \label{1eqInjective2}
    \text{if} \quad S \in E' \quad \text{or} \quad T \in E' \quad
    \text{then} \quad g'_S = g'_T \quad \text{implies} \quad S = T.
  \end{equation}
  \begin{claim} \label{c1}
    For all $S,T \in E$ we have:
    \[
    S \cap T \ne \emptyset \quad \Longleftrightarrow \quad \phi(S)
    \cap \phi(T) \ne \emptyset.
    \]
  \end{claim}
  The claim is clearly true if $S = T$, so we may assume that $S \ne
  T$. To prove the implication ``$\Rightarrow$'', write $S = \{i,j\}$,
  $T = \{i,k\}$ with $i,j,k \in \N$ pairwise distinct. The (multi-)
  set $\{(g\ss{i,j},w_k),(g\ss{i,k},w_j),(g\ss{j,k},w_i)\}$ occurs in
  ${\mathcal T}_G$ and therefore also in ${\mathcal T}_{G'}$. Hence
  there exist pairwise distinct $r,s,t \in \N$ with
  \[
  \{g\ss{i,j},g\ss{i,k},g\ss{j,k}\} =
  \{g'\ss{r,s},g'\ss{r,t},g'\ss{s,t}\}.
  \]
  By~\eqref{1eqPhi} this implies
  \[
  \{g'\sps{i,j},g'\sps{i,k}\} \subseteq
  \{g'\ss{r,s},g'\ss{r,t},g'\ss{s,t}\},
  \]
  so by~\eqref{1eqInjective2}
  \[
  \ps{i,j},\ps{i,k} \in \{\{r,s\},\{r,t\},\{s,t\}\}.
  \]
  This implies $\ps{i,j} \cap \ps{i,k} \ne \emptyset$.

  To prove the converse, write $\phi(S) = \{i,j\}$ and $\phi(T) =
  \{i,k\}$ with $i,j,k \in \N$ pairwise distinct. Since
  $\{(g'\ss{i,j},w'_k),(g'\ss{i,k},w'_j),(g'\ss{j,k},w'_i)\}$ occurs
  in ${\mathcal T}_{G'}$ there exist pairwise distinct $r,s,t \in \N$
  such that
  \[
  \{g'\ss{i,j},g'\ss{i,k},g'\ss{j,k}\} =
  \{g\ss{r,s},g\ss{r,t},g\ss{s,t}\}.
  \]
  But $g_S = g'\ss{i,j}$ and $g_T = g'\ss{i,k}$ are elements in the
  set on the left hand side, so by~\eqref{1eqInjective} we see that $S
  \cap T \ne \emptyset$. This completes the proof of \clref{c1}.

  For $i \in \N$, write
  \[
  N_i := \{j \in \N \mid \{i,j\} \in E\}.
  \]
  Define
  \[
  I := \left\{i \in \N \mid |N_i| \ge 2\right\}.
  \]
  \begin{claim} \label{c2}
    For every $i \in I$, the intersection $\bigcap_{j \in N_i}
    \ps{i,j}$ has precisely one element.
  \end{claim}
  Indeed, the injectivity of~$\phi$ implies that the intersection has
  at most one element. Let $i \in I$ and choose $j,k \in N_i$
  distinct. By \clref{c1}, $\ps{i,j}$ and $\ps{i,k}$ have non-empty
  intersection, so we can write $\ps{i,j} = \{r,s\}$ and $\ps{i,k} =
  \{r,t\}$ with $r,s,t \in \N$ pairwise distinct. We have to show that
  for every $l \in N_i$ we have $r \in \ps{i,l}$. Assume the contrary.
  Since $\ps{i,l}$ has a non-empty intersection with $\{r,s\}$ and
  with $\{r,t\}$ by \clref{c1}, $r \notin \ps{i,l}$ implies $\ps{i,l}
  = \{s,t\}$. Since
  $\{(g'\ss{r,s},w'_t),(g'\ss{r,t},w'_s),(g'\ss{s,t},w'_r)\}$ occurs
  in ${\mathcal T}_{G'}$, there exist pairwise distinct $u,v,w \in \N$
  such that
  \[
  \{g\ss{u,v},g\ss{u,w},g\ss{v,w}\} =
  \{g'\ss{r,s},g'\ss{r,t},g'\ss{s,t}\} =
  \{g\ss{i,j},g\ss{i,k},g\ss{i,l}\},
  \]
  hence by~\eqref{1eqInjective} we have
  $\left\{\{u,v\},\{u,w\},\{v,w\}\right\} =
  \left\{\{i,j\},\{i,k\},\{i,l\}\right\}$. Forming the intersection of
  both sides yields the contradiction $\emptyset = \{i\}$. This proves
  \clref{c2}.

  By \clref{c2}, we may define a map $\map{\pi}{I}{\N}$ by
  \begin{equation} \label{1eqPi}
    \bigcap_{j \in N_i} \ps{i,j} = \{\pi(i)\} \quad \text{for all}
    \quad i \in I.
  \end{equation}
  \begin{claim} \label{c3}
    The map~$\pi$ is injective.
  \end{claim}
  To prove this claim, assume that there exist $i,r \in I$ with $i \ne
  r$ such that $\pi(i) = \pi(r)$. Choose $j,k \in N_i$ distinct and
  $s,t \in N_r$ distinct. By~\eqref{1eqPi} we have
  \begin{equation} \label{1eqIntersect}
    \pi(i) = \pi(r) \in \ps{i,j} \cap \ps{i,k} \cap \ps{r,s} \cap
    \ps{r,t},
  \end{equation}
  so
  \[
  \{i,j\} \cap \{r,s\} \ne \emptyset \quad \text{and} \quad \{i,j\}
  \cap \{r,t\} \ne \emptyset
  \]
  by \clref{c1}. This implies $i = s$ or $i = t$ or $j = r$. In all
  cases, $\{i,r\} \in E$ follows, thus we may specify our choice
  of~$k$ and~$t$ further by setting $k = r$ and $t = i$. This choice
  implies $j \ne r$ and $s \ne i$. Now $\{i,j\} \cap \{r,s\} \ne
  \emptyset$ implies $s = j$. From~\eqref{1eqIntersect} we obtain
  \begin{equation} \label{1eqIntersect2}
    \pi(i) \in \ps{i,j} \cap \ps{i,r} \cap \ps{j,r}.
  \end{equation}
  Since ${\mathcal T}_{G'} = {\mathcal T}_G$, there exist pairwise
  distinct $u,v,w \in \N$ with $\{g\ss{i,j},g\ss{i,r},g\ss{j,r}\} =
  \{g'\ss{u,v},g'\ss{u,w},g'\ss{v,w}\}$. Using~\eqref{1eqPhi} we obtain
  $\{g'\sps{i,j},g'\sps{i,r},g'\sps{j,r}\} =
  \{g'\ss{u,v},g'\ss{u,w},g'\ss{v,w}\}$, so
  \[
  \left\{\ps{i,j},\ps{i,r},\ps{j,r}\right\} =
  \left\{\{u,v\},\{u,w\},\{v,w\}\right\}
  \]
  by~\eqref{1eqInjective2}.  Forming the intersection of both sides
  yields a contradiction with~\eqref{1eqIntersect2}. This proves
  \clref{c3}.
  \begin{claim} \label{c4}
    For all $i,j \in I$ with $\{i,j\} \in E$ we have
    \[
    \ps{i,j} = \{\pi(i),\pi(j)\}.
    \]
  \end{claim}
  Indeed, $\{i,j\} \in E$ implies $j \in N_i$ and $i \in N_j$, hence
  by~\eqref{1eqPi}
  \[
  \pi(i) \in \ps{i,j} \quad \text{and} \quad \pi(j) \in \ps{i,j}.
  \]
  Since~$i$ and~$j$ are distinct, \clref{c3} implies $\pi(i) \ne
  \pi(j)$. Now \clref{c4} follows, since $\ps{i,j}$ has precisely two
  elements.
  \begin{claim} \label{c4b}
    For every $i \in I$ we have
    \[
    w'_{\pi(i)} = w_i.
    \]
    (Recall that $w_i$ is the weight of node~$i$, see before the
    statement of the theorem.)
  \end{claim}
  For the proof, let $i \in I$ and choose $j,k \in N_i$ distinct.
  Because ${\mathcal T}_{G'} = {\mathcal T}_G$, there exist pairwise
  distinct $r,s,t \in \N$ such that
  \begin{equation} \label{1eq4b}
    \{(g\ss{i,j},w_k),(g\ss{i,k},w_j),(g\ss{j,k},w_i)\} =
    \{(g'\ss{r,s},w'_t),(g'\ss{r,t},w'_s),(g'\ss{s,t},w'_r)\}.
  \end{equation}
  By~\eqref{1eqPhi} it follows that
  \[
  \{g'\sps{i,j},g'\sps{i,k}\} \subseteq
  \{g'\ss{r,s},g'\ss{r,t},g'\ss{s,t}\}.
  \]
  By~\eqref{1eqPi} we can write $\phi(\{i,j\}) = \{\pi(i),u\}$ and
  $\phi(\{i,k\}) = \{\pi(i),v\}$ with $u,v \in \N \setminus
  \{\pi(i)\}$ distinct. Now the above inclusion together
  with~\eqref{1eqInjective2} implies
  \[
  \left\{\{\pi(i),u\},\{\pi(i),v\}\right\} \subseteq
  \left\{\{r,s\},\{r,t\},\{s,t\}\right\},
  \]
  so $\{\pi(i),u,v\} = \{r,s,t\}$. With this, \eqref{1eq4b} becomes
  \[
  \{(g\ss{i,j},w_k),(g\ss{i,k},w_j),(g\ss{j,k},w_i)\} =
  \{(g'\ss{\pi(i),u},w'_v),(g'\ss{\pi(i),v},w'_u),(g'\ss{u,v},w'_{\pi(i)})\}.
  \]
  We have $g\ss{i,j} = g'\ss{\pi(i),u}$ and $g\ss{i,k} =
  g'\ss{\pi(i),v}$. Since $g\ss{i,j}$, $g\ss{i,k}$, $g\ss{j,k}$ are
  pairwise distinct, it follows that $(g\ss{j,k},w_i) =
  (g'\ss{u,v},w'_{\pi(i)})$. Hence $w_i = w'_{\pi(i)}$, proving
  \clref{c4b}.

  Set $B := \{j \in \N \mid |N_j| = 1\}$ and
  \[
  B_0 := \left\{j \in B \mid N_j \subseteq I\right\}, \quad B_1 :=
  \left\{j \in B \mid N_j \subseteq \N \setminus I\right\},
  \]
  so $B$ is the disjoint union of $B_0$ and $B_1$. For every $j \in
  B_0$ there exists a unique $i \in N_j$, and this~$i$ lies in $I$. We
  also have $j \in N_i$, so $\pi(i) \in \ps{i,j}$ by~\eqref{1eqPi}. We
  can thus extend~$\pi$ to a map $\map{\pi}{I \cup B_0}{\N}$ by
  setting $\ps{i,j} = \{\pi(i),\pi(j)\}$.
  \begin{claim} \label{c5}
    The map $\map{\pi}{I \cup B_0}{\N}$ is injective, and for $i,j \in
    I \cup B_0$ with $\{i,j\} \in E$ we have
    \[
    \ps{i,j} = \{\pi(i),\pi(j)\}.
    \]
  \end{claim}
  To prove this claim, take $i,j \in I \cup B_0$ distinct. If $i,j \in
  I$, then both assertions follow from Claims~\ref{c3} and~\ref{c4}.

  Next, assume $i \in I$ and $j \in B_0$. The second assertion follows
  from the definition of~$\pi$. By way of contradiction, assume that
  $\pi(i) = \pi(j)$. Choose $k,l \in N_i$ distinct and write $N_j =
  \{r\}$. Then
  \[
  \pi(i) = \pi(j) \in \ps{i,k} \cap \ps{i,l} \cap \ps{j,r},
  \]
  so \clref{c1} implies $i = r$ or $\{k,l\} = \{j,r\}$. The second
  case implies $j \in N_i$, or, equivalently, $i \in N_j$, so $i = r$
  in this case, too. We obtain
  \[
  \ps{j,r} = \{\pi(j),\pi(r)\} = \{\pi(j),\pi(i)\},
  \]
  which implies $\pi(i) \ne \pi(j)$ since $\ps{j,r}$ has two elements.

  Finally, assume $i,j \in B_0$ distinct. Then $\{i,j\} \notin E$ by
  the definition of $B_0$, so the second assertion of \clref{c5} holds
  automatically. Write $N_i = \{k\}$ and $N_j = \{l\}$. Then
  \[
  \ps{i,k} = \{\pi(i),\pi(k)\} \quad \text{and} \quad \ps{j,l} =
  \{\pi(j),\pi(l)\}
  \]
  by the definition of~$\pi$, hence $\pi(i) = \pi(j)$ would imply
  $\{i,k\} \cap \{j,l\} \ne \emptyset$ by \clref{c1}, so $k = l$ since
  $i,j \in B_0$ but $k,l \in I$. But then
  \[
  \ps{i,k} = \{\pi(i),\pi(k)\} = \{\pi(j),\pi(l)\} = \ps{j,l} =
  \ps{j,k},
  \]
  contradicting the injectivity of~$\phi$. Therefore $\pi(i) \ne
  \pi(j)$, which completes the proof of \clref{c5}.
  \begin{claim} \label{c5b}
    For every $i \in I \cup B_0$ we have $w_i = w'_{\pi(i)}$.
  \end{claim}
  The proof is similar to the one of \clref{c4b}. First, if $i \in I$
  this is \clref{c4b}. So assume $i \in B_0$ and write $N_i =
  \{j\}$. Since $j \in I$, we can choose $k \in N_j \setminus \{i\}$.
  Because ${\mathcal T}_{G'} = {\mathcal T}_G$, there exist pairwise
  distinct $r,s,t \in \N$ such that
  \begin{equation} \label{1eq5b}
    \{(g\ss{i,j},w_k),(g\ss{i,k},w_j),(g\ss{j,k},w_i)\} =
    \{(g'\ss{r,s},w'_t),(g'\ss{r,t},w'_s),(g'\ss{s,t},w'_r)\}.
  \end{equation}
  By~\eqref{1eqPhi}, it follows that
  \[
  \{g'\sps{i,j},g'\sps{j,k}\} \subseteq
  \{g'\ss{r,s},g'\ss{r,t},g'\ss{s,t}\}.
  \]
  By \clref{c5} we have $\phi(\{i,j\}) = \{\pi(i),\pi(j)\}$ and
  $\phi(\{j,k\}) = \{\pi(j),\pi(k)\}$. Now the above inclusion
  together with~\eqref{1eqInjective2} implies
  \[
  \left\{\{\pi(i),\pi(j)\},\{\pi(j),\pi(k)\}\right\} \subseteq
  \left\{\{r,s\},\{r,t\},\{s,t\}\right\},
  \]
  so $\{\pi(i),\pi(j),\pi(k)\} = \{r,s,t\}$. With this, \eqref{1eq5b}
  becomes
  \begin{multline*}
    \{(g\ss{i,j},w_k),(g\ss{i,k},w_j),(g\ss{j,k},w_i)\} = \\
    \{(g'\ss{\pi(i),\pi(j)},w'_{\pi(k)}),
    (g'\ss{\pi(i),\pi(k)},w'_{\pi(j)}),
    (g'\ss{\pi(j),\pi(k)},w'_{\pi(i)})\}.
  \end{multline*}
  We have $g\ss{i,j} = g'\ss{\pi(i),\pi(j)}$ and $g\ss{j,k} =
  g'\ss{\pi(j),\pi(k)}$. Since $g\ss{i,j}$, $g\ss{i,k}$, $g\ss{j,k}$
  are pairwise distinct, it follows that $(g\ss{j,k},w_i) =
  (g'\ss{\pi(j),\pi(k)},w'_{\pi(i)})$. Hence $w_i = w'_{\pi(i)}$,
  proving \clref{c5b}.

  Now we turn our attention to $B_1$.
  \begin{claim} \label{c5c}
    Let $i,j \in B_1$ with $\{i,j\} \in E$, and write $\phi(\{i,j\}) =
    \{l,m\}$. Then
    \[
    \{w_i,w_j\} = \{w'_l,w'_m\}.
    \]
  \end{claim}
  For the proof, choose $k \in \N \setminus \{i,j\}$ arbitrary.
  Because ${\mathcal T}_{G'} = {\mathcal T}_G$, there exist pairwise
  distinct $r,s,t \in \N$ such that~\eqref{1eq5b} holds. We have
  $g'\ss{l,m} = g\ss{i,j}$, so
%  But $i,j \in
%  B_1$ implies $g\ss{i,k} = g\ss{j,k} = 0$, and we have $g\ss{i,j} =
%  g'\ss{l,m}$, so~\eqref{1eq5b} becomes
%  \begin{equation} \label{1eq5c}
%    \{(g'\ss{l,m},w_k),(0,w_j),(0,w_i)\} =
%    \{(g'\ss{r,s},w'_t),(g'\ss{r,t},w'_s),(g'\ss{s,t},w'_r)\}.
%  \end{equation}
  by~\eqref{1eqInjective2} we conclude $\{l,m\} \in
  \left\{\{r,s\},\{r,t\},\{s,t\}\right\}$, and by interchanging the
  roles of~$r$, $s$ and~$t$ we may assume that $\{l,m\} = \{r,s\}$. So
  we must have $(g'\ss{l,m},w_k) = (g'\ss{r,s},w'_t)$,
  and~\eqref{1eq5b} becomes
  \[
  \{(g\ss{i,k},w_j),(g\ss{j,k},w_i)\} =
  \{(g'\ss{l,t},w'_m),(g'\ss{m,t},w'_l)\}.
  \]
  \clref{c5c} follows from this.

  For each $i \in B_1$ we have $N_i = \{j\}$ with $j \notin I$, so $j
  \in B_1$ and thus $N_j = \{i\}$. Thus we may extend~$\pi$ to a map
  $\map{\pi}{I \cup B}{\N}$ such that for all $i,j \in B_1$ with
  $\{i,j\} \in E$ we have $\ps{i,j} =
  \{\pi(i),\pi(j)\}$. This condition still holds if the values of
  $\pi(i)$ and $\pi(j)$ are swapped. By~\clref{c5c} it is possible to
  specify the extension of~$\pi$ further such that
  \[
  w_i = w'_{\pi(i)}
  \]
  for all $i \in B_1$.

  \begin{claim} \label{c6}
    The map $\map{\pi}{I \cup B}{\N}$ is injective, and for $i,j \in
    I \cup B$ with $\{i,j\} \in E$ we have
    \[
    \ps{i,j} = \{\pi(i),\pi(j)\}.
    \]
  \end{claim}
  For the proof, take $i,j \in I \cup B$ distinct. If $i,j \in I \cup
  B_0$, then both assertions follow from \clref{c5}.

  Next, assume $i \in I \cup B_0$ and $j \in B_1$. Then $\{i,j\} \in
  E$ is impossible, so the second assertion is automatic. It follows
  from the definitions of $I$, $B_0$ and $B_1$ that there exist $k \in
  I \cup B_0$ and $l \in B_1$ such that $\{i,k\}, \{j,l\} \in E$. From
  \clref{c5} and the definition of~$\pi$ we have
  \[
  \ps{i,k} = \{\pi(i),\pi(k)\} \quad \text{and} \quad \ps{j,l} =
  \{\pi(j),\pi(l)\}.
  \]
  But $\{i,k\} \cap \{j,l\} \subseteq (I \cup B_0) \cap B_1 =
  \emptyset$, so
  \[
  \{\pi(i),\pi(k)\} \cap \{\pi(j),\pi(l)\} = \ps{i,k} \cap \ps{j,l} =
  \emptyset
  \]
  by \clref{c1}. This implies $\pi(i) \ne \pi(j)$.

  Finally, assume $i,j \in B_1$. If $\{i,j\} \in E$, then $\ps{i,j} =
  \{\pi(i),\pi(j)\}$ by the definition of~$\pi$. In that case, $\pi(i)
  \ne \pi(j)$ since $\ps{i,j}$ has two elements. On the other hand, if
  $\{i,j\} \notin E$, there exist $k,l \in B_1$ with $\{i,k\},\{j,l\}
  \in E$. The definition of $B$ implies $k \ne l$, and $\{i,j\} \notin
  E$ implies $i \ne l$ and $k \ne j$. Hence $\{i,k\} \cap \{j,l\} =
  \emptyset$. By \clref{c1}, this implies
  \[
  \{\pi(i),\pi(k)\} \cap \{\pi(j),\pi(l)\} = \ps{i,k} \cap \ps{j,l} =
  \emptyset,
  \]
  so $\pi(i) \ne \pi(j)$. This completes the proof of \clref{c6}.

  In summary, we have an injective map $\map{\pi}{I \cup B}{\N}$ with
  the property that $w_i = w'_{\pi(i)}$ for all $i \in I \cup B$. Form
  the union (with adding multiplicities) of all multisets ${\mathfrak
  t}\ss{i,j,k}$ lying in ${\mathcal T}_G$, and then take the multiset
  consisting of the second components of all $(g\ss{i,j},w_k)$ lying
  in this union. This yields the multiset consisting of all $w_i$, $i
  \in \N$, counted $\binom{n-1}{2}$ times for each $i \in \N$. Since
  ${\mathcal T}_{G'} = {\mathcal T}_G$, this implies that the multiset
  of all $w_i$, $i \in \N$ and the multiset of all $w'_j$, $j \in \N$,
  coincide. Therefore we can extend~$\pi$ to obtain a bijection
  $\map{\pi}{\N}{\N}$ such that
  \begin{equation} \label{1eqW}
    w_i = w'_{\pi(i)}
  \end{equation}
  holds for all $i \in \N$. This map~$\pi$ induces a bijection
  $\map{\phi_\pi}{P}{P}$ defined by
  \[
  \phi_\pi\left(\{i,j\}\right) = \{\pi(i),\pi(j)\} \quad \text{for}
  \quad \{i,j\} \in P.
  \]
  For $\{i,j\} \in E$ we have $i,j \in I \cup B$ by the definition of
  $I$ and $B$, hence \clref{c6} says that the restriction of
  $\phi_\pi$ to $E$ is~$\phi$. Thus for $\{i,j\} \in E$ we have
  \[
  g\ss{i,j} = g'\sps{i,j} = g'_{\phi_\pi\left(\{i,j\}\right)} =
  g'\ss{\pi(i),\pi(j)},
  \]
  where the first equation follows from~\eqref{1eqPhi}. Finally, take
  $i,j,k \in \N$ pairwise distinct with $\{i,k\}, \{j,k\} \in E$. By
  the above, this implies
  \begin{equation} \label{1eqC}
    g\ss{i,k} = g'\ss{\pi(i)\pi(k)} \quad \text{and} \quad g\ss{j,k} =
    g'\ss{\pi(j)\pi(k)}.
  \end{equation}
  Since ${\mathcal T}_{G'} = {\mathcal T}_G$,
  there exist pairwise distinct $r,s,t \in \N$ such that
  \begin{equation} \label{1eqC2}
    \{g\ss{i,j},g\ss{i,k},g\ss{j,k}\} =
    \{g'\ss{r,s},g'\ss{r,t},g'\ss{s,t}\}.
  \end{equation}
  Using~\eqref{1eqC} and~\eqref{1eqInjective}, we obtain
  \[
  \left\{\{\pi(i),\pi(k)\},\{\pi(j),\pi(k)\}\right\} \subseteq
  \left\{\{r,s\},\{r,t\},\{s,t\}\right\},
  \]
  so $\{r,s,t\} = \{\pi(i),\pi(j),\pi(k)\}$, and~\eqref{1eqC2} becomes
  \[
  \{g\ss{i,j},g\ss{i,k},g\ss{j,k}\} =
  \{g'\ss{\pi(i),\pi(j)},g'\ss{\pi(i),\pi(k)},g'\ss{\pi(j),\pi(k)}\}.
  \]
  Since the set on the left side has three distinct elements, we
  conclude, using~\eqref{1eqC}, that $g\ss{i,j} =
  g'\ss{\pi(i),\pi(j)}$. This completes the proof.
%  On the other hand, for $\{i,j\} \in P \setminus E$ we have
%  $\{\pi(i),\pi(j)\} = \phi_\pi\left(\{i,j\}\right) \notin E'$ (since
%  $\phi_\pi$ maps $E$ onto $E'$), so
%  \[
%  g'\ss{\pi(i),\pi(j)} = 0 = g\ss{i,j}.
%  \]
%  So for all $\{i,j\} \in P$ we have $g\ss{i,j} =
%  g'\ss{\pi(i),\pi(j)}$. Together with~\eqref{1eqW} this means that
%  $G'$ and $G$ are isomorphic, completing the proof.
\end{proof}

\subsection{Distinct weights} \label{1sDistinct}

We say that a graph $G$ with weighted edges is \df{reconstructible
  from the distribution of subtriangles} if $G$ is reconstructible
from $\mathcal T$. In other words, we demand that every graph $G'$
with the same distribution of subtriangles as $G$ is isomorphic to
$G$. 

\begin{ex} \label{1exReconstrct}
  Figure~\ref{1fCounter} on page~\pageref{1fCounter} shows a pair of
  graphs which are not isomorphic, but have the same distribution of
  subtriangles. Each edge that is drawn represents an edge of
  weight~1, and an edge which is not drawn represents weight~0. The
  node weights can all be taken to be~0. Simple counting reveals that
  in both graphs there are two subtriangles with all edge weights~0,
  four subtriangles with one non-zero edge weight, four subtriangles
  with two non-zero edge weights, and no subtriangle with all weights
  non-zero. So there exist graphs which are not reconstructible from
  the distribution of subtriangles. In fact, Figure~\ref{1fCounter}
  gives the simplest such example.
\end{ex}

Let $0 \in X$ be some distinguished weight, so $g\ss{i,j} = 0$ may be
interpreted as saying that the nodes~$i$ and~$j$ are not connected.
The hypotheses of the following theorem sound a bit technical. For
that reason we formulate a special case as \cref{1cReconstruct}, where
the hypotheses are easier to state (and to remember), so readers might
wish to read \cref{1cReconstruct} first.  \exref{1exReconstruct} is a
typical example where the hypotheses of \tref{1tReconstruct}, but not
those of \cref{1cReconstruct}, are satisfied.

\begin{theorem} \label{1tReconstruct}
  Let $G$ be a graph with $n \ge 3$ nodes, with edge weights
  $g\ss{i,j}$ and node weights $w_i = g\ss{i,i}$. Write
  \[
  P := \left\{\{i,j\} \mid 1 \le i < j \le n\right\}
  \]
  and
  \[
  E := \left\{S \in P \mid g_T \ne g_S \ \text{for all} \ T \in P
    \setminus \{S\}\right\}.
  \]
  Assume that for every $\{i,j\} \in P$ at least one of the following
  conditions holds:
  \begin{enumerate}
    \renewcommand{\theenumi}{\roman{enumi}}
  \item $\{i,j\} \in E$,
  \item there exists $k \in \N \setminus \{i,j\}$ such that $\{i,k\}
    \in E$ and $\{j,k\} \in E$, or
  \item $g\ss{i,j} = 0$.
  \end{enumerate}
  Then $G$ is reconstructible from the distribution of subtriangles
\end{theorem}

\begin{proof}
  Let $G'$ be a graph with ${\mathcal T}_{G'} = {\mathcal T}_G$, and
  let $\pi \in S_n$ be a bijection as given by \lref{1lReconstruct}.
  Writing $\map{\phi_\pi}{P}{P}$ for the map induced by~$\pi$, we
  obtain
  \[
  g\ss{i,j} = g'_{\phi_\pi(\{i,j\})} = g'\ss{\pi(i),\pi(j)}
  \]
  for all $\{i,j\} \in P$ satisfying condition~(i) or~(ii) of the
  theorem. In particular, if $g_S$ is non-zero for an $S \in P$, then
  the same is true for $g'_{\phi_\pi(S)}$. But since the multisets of
  all $g_S$ and of all $g'_S$ coincide (see the beginning of the proof
  of \lref{1lReconstruct}), it follows that if $g_S$ is zero, then
  $g'_{\phi_\pi(S)}$ is zero, too. So for $\{i,j\} \in P$ satisfying
  condition~(iii) we have
  \[
  g\ss{i,j} = 0 = g'_{\phi_\pi(\{i,j\})} = g'\ss{\pi(i),\pi(j)}.
  \]
  In summary, we have $g\ss{i,j} = g'\ss{\pi(i),\pi(j)}$ for all
  $\{i,j\} \in P$ and $w_i = w'_{\pi(i)}$ for all $i \in \N$, which
  means that $G$ and $G'$ are isomorphic.
\end{proof}

The following corollary is just a weaker version of
\tref{1tReconstruct}. The hypothesis means that~0 is the only weight
between distinct edges that may occur repeatedly, or, by interpreting
weight 0 as ``not connected'', that the edge weights between pairs of
connected nodes are pairwise distinct.

\begin{cor} \label{1cReconstruct}
  Let $G$ be a graph with $n \ge 3$ nodes, with edge weights
  $g\ss{i,j}$ and node weights $w_i = g\ss{i,j}$. Assume that for
  $i,j,k,l \in \N$ with $i \ne j$ and $k \ne l$ we have that
  $g\ss{i,j} = g\ss{k,l} \ne 0$ implies $\{i,j\} = \{k,l\}$. Then $G$
  is reconstructible from the distribution of subtriangles.
\end{cor}

\begin{ex} \label{1exReconstruct}
  Figure~\ref{1fReconstruct} shows a graph with~5 nodes to which
  \tref{1tReconstruct} is applicable, but \cref{1cReconstruct} is not.
  Here~$a$, $b$, $c$ and~$d$ denote pairwise distinct, non-zero edge
  weights, and edges which are not drawn are to be understood as
  having weight~0. There are no node weights. (In fact, this example
  would also be valid with node weights assigned arbitrarily.)
\end{ex}

\providecommand{\point}{\circle*{2.5}}
\newcommand{\attach}[1]{\small $#1$}
\Figure{
  \fbox{
    \begin{picture}(120,120)
      \put(60,20){\point}
      \put(20,60){\point}
      \put(60,60){\point}
      \put(100,60){\point}
      \put(60,100){\point}
      \put(23,63){\line(1,1){34}}
      \put(33,80){\attach{$d$}}
      \put(23,60){\line(1,0){34}}
      \put(38,54){\attach{$a$}}
      \put(63,60){\line(1,0){34}}
      \put(78,54){\attach{$c$}}
      \put(60,63){\line(0,1){34}}
      \put(61,72){\attach{$b$}}
      \put(63,97){\line(1,-1){34}}
      \put(81,80){\attach{$d$}}
%      \put(60,23){\line(0,1){34}}
%      \put(61,38){\attach{$0$}}
%      \put(63,23){\line(1,1){34}}
%      \put(81,35){\attach{$0$}}
%      \put(23,57){\line(1,-1){34}}
%      \put(33,35){\attach{$0$}}
    \end{picture}
  }
%  \fbox{
%    \begin{picture}(120,120)
%      \put(60,20){\point}
%      \put(20,60){\point}
%      \put(60,60){\point}
%      \put(100,60){\point}
%      \put(60,100){\point}
%      \put(20,60){\line(1,1){40}}
%      \put(33,80){\attach{$d$}}
%      \put(20,60){\line(1,0){40}}
%      \put(38,54){\attach{$a$}}
%      \put(60,60){\line(1,0){40}}
%      \put(78,54){\attach{$c$}}
%      \put(60,60){\line(0,1){40}}
%      \put(61,72){\attach{$b$}}
%      \put(60,100){\line(1,-1){40}}
%      \put(81,80){\attach{$d$}}
%    \end{picture}
%  }
}
{A graph to which \tref{1tReconstruct} is applicable (edges which are
  not drawn have weight~0)}
{1fReconstruct}

We conclude this section by asking whether distributions of
subtriangles can be represented in a practical, computer friendly way.
Multisets are not very practical, since the are hard to visualize and
to compare.  Whether a better representation exists depends on the set
$X$ in which the weights $g\ss{i,j}$ and $w_i$ lie. If $X$ is finite,
we may assume $X = \{1 \upto r-1\}$ with an integer~$r$. A pair
$(g\ss{i,j},w_k)$ can then be uniquely represented by the single
integer $g\ss{i,j} r + w_k < r^2$. Instead of representing a
subtriangle ${\mathfrak t}\ss{i,j,k}$ as a multiset, we can order the
elements (represented as integers less than $r^2$) by size, so we get
$0 \le a_1 \le a_2 \le a_3 < r^2$. These can be uniquely represented
by $a_1 r^4 + a_2 r^2 + a_3 < r^6$. Continuing this way, ${\mathcal
  T}_G$ may be uniquely represented as a single integer between~0 and
$r^{6 \binom{n}{3}} - 1$.

However, the situation becomes much more tricky if $X$ is
infinite. Even for $X = \RR$, it is far from clear that there is a
good way of representing the distribution of subtriangles, and for $X
= \RR^d$ it becomes harder still. That is why we turn to simpler,
one-dimensional distributions in \sref{3s:Mimi}.

\subsection{Some statistics} \label{1sStats}

The hypotheses of \tref{1tReconstruct} and \cref{1cReconstruct} are,
very roughly speaking, that the edge weights are sufficiently
distinct. It is therefore clear that these hypothesis will tend to be
met if the weights take a large range of values, and vice versa.  But
even if the weights lie in a small set, we can hope that graphs are
reconstructible from the distribution of subtriangles even if our
theorem fails to guarantee that. To get some idea of to how many graphs
this applies, we ran some computer experiments. The results are given
in Table~\ref{1taStats}.

For several values of~$n$ and~$m$, we considered all graphs with~$n$
points, where the edge weights appearing in the graph are precisely
all integers between~0 and $m-1$. We considered graphs without node
weights. The third column in the table gives the number of graphs
which are not reconstructible from the distribution of subtriangles,
divided by the number of all graphs. For example, 98.8\% of all graphs
with~5 nodes and~6 different edge weights are reconstructible from the
distribution of subtriangles. Even if most graphs are not
reconstructible (as in the case of~6 nodes and~2 different weights),
it is still possible that there only exists a small set of pairs of
graphs which are not isomorphic but have the same distribution of
subtriangles. In fact, this tends to be the case. For example, there
are $2^{15} = 32768$ simple graphs of 6 nodes, affording $2^{30}
\approx 10^9$ ordered pairs of graphs. Among those, we found precisely
7680960 pairs of two graphs which are non-isomorphic but have the same
distribution of triangles. We interpret this by saying that for simple
graphs with~6 nodes, the error probability of testing by subtriangles
is $7680960/2^{30} \approx 7.15 \cdot 10^{-3}$. These error
probabilities are given in the fourth column of Table~\ref{1taStats}.

We compared this to the possibilities of discriminating non-isomorphic
graphs by another invariant, the spectrum. This is perhaps the
best-known graph invariant. By definition, the spectrum of a graph is
the set of eigenvalues, with multiplicities, of the adjacency matrix.
The fifth column of Table~\ref{1taStats} contains the error
probabilities when one tries to discriminate non-isomorphic graphs by
using the spectrum. So the entry is the number of all pairs of
cospectral graphs with a given number of nodes and weights, divided by
the number of all pairs. But the spectrum and the distribution of
subtriangles should not be regarded as competing invariants, because
together they work best, as the last column of the table shows. This
column gives the error probabilities one gets when combining the
spectrum and the distribution of subtriangles.

The computation were all done by using the computer algebra system Magma
(see \mycite{magma}).

\begin{table}[h]
  \centering
  \newcommand{\err}{P_{\operatorname{err}}}
  \begin{tabular}{|c|c||c||c|c|c|}
    \hline
    nodes & weights & non-reconst. & $\err$ triangles & $\err$ spectrum & $\err$
    combination \\
    \hline
     4 & 2 &
    0
    &
    0
    &
    0
    &
    0
    \\
    \hline
     4 & 3 &
    0
    &
    0
    &
    0
    &
    0
    \\
    \hline
     4 & 4 &
    0
    &
    0
    &
    0
    &
    0
    \\
    \hline
     4 & 5 &
    0
    &
    0
    &
    3.56 $\cdot 10^{-4}$
    &
    0
    \\
    \hline
     4 & 6 &
    0
    &
    0
    &
    0
    &
    0
    \\
    \hline
     5 & 2 &
    0.15
    &
    3.45 $\cdot 10^{-3}$
    &
    1.44 $\cdot 10^{-4}$
    &
    0
    \\
    \hline
     5 & 3 &
    0.15
    &
    2.83 $\cdot 10^{-4}$
    &
    3.85 $\cdot 10^{-5}$
    &
    2.30 $\cdot 10^{-6}$
    \\
    \hline
     5 & 4 &
    0.078
    &
    1.28 $\cdot 10^{-5}$
    &
    3.20 $\cdot 10^{-6}$
    &
    0
    \\
    \hline
     5 & 5 &
    0.034
    &
    8.55 $\cdot 10^{-7}$
    &
    4.91 $\cdot 10^{-7}$
    &
    2.21 $\cdot 10^{-9}$
    \\
    \hline
     5 & 6 &
    0.012
    &
    8.64 $\cdot 10^{-8}$
    &
    1.07 $\cdot 10^{-7}$
    &
    0
    \\
    \hline
     5 & 7 &
    0.0026
    &
    1.03 $\cdot 10^{-8}$
    &
    3.07 $\cdot 10^{-8}$
    &
    0
    \\
    \hline
     6 & 2 &
    0.63
    &
    7.15 $\cdot 10^{-3}$
    &
    9.56 $\cdot 10^{-5}$
    &
    4.02 $\cdot 10^{-5}$
    \\
    \hline
     6 & 3 &
    0.62
    &
    6.68 $\cdot 10^{-5}$
    &
    4.28 $\cdot 10^{-6}$
    &
    8.12 $\cdot 10^{-7}$
    \\
    \hline
     7 & 2 &
    0.93
    &
    6.55 $\cdot 10^{-3}$
    &
    9.36 $\cdot 10^{-5}$
    &
    6.55 $\cdot 10^{-5}$
    \\
    \hline
     8 & 2 &
    0.99
    &
    3.97 $\cdot 10^{-3}$
    &
    1.58 $\cdot 10^{-5}$
    &
    1.40 $\cdot 10^{-5}$
    \\
    \hline
  \end{tabular}
  \caption{Ratio of graphs which are not reconstructible from
    subtriangles, and error probabilities using subtriangles and
    spectra}
  \label{1taStats}
\end{table}

%%% Local Variables: 
%%% mode: latex
%%% TeX-master: "graphs"
%%% End: 

\section{Reconstructibility from one-dimensional distributions} 
\label{3s:Mimi}

In this section, we concentrate on the case where the weights of a graph take on real (vector) values.
The representations we construct are in terms of simple one-dimensional distributions.
As we show in the following, except for a set of measure zero of graphs,
the representations we propose are lossless.

Most of the claims in this section are based on the following lemma.
A proof can be found in \cite{Boutin.Kemper03}.
\begin{lemma}
 \label{L:permutations}
 Let $n\geq 5$.
 Consider the action of the permutation group 
 $S_{\binom{n}{2}}$ on the set of pairs
 $\{ \{ i,j\} | i,j=1,\ldots, n , i\neq j \}$.
 Let $\varphi\in S_{\binom{n}{2}}$.
 Then there exists a permutation $\pi\in S_n$
 such that
 $\varphi \left( \{ i,j\}\right) =\{ \pi(i),\pi(j)\}$, for every $i\neq j$,
 if and only if
 $\varphi (\{i,j \})\cap \varphi (\{j,k \})\neq \emptyset$
 for every pairwise distinct $i,j,k\in \{1,\ldots,n \}$. 
\end{lemma}

%%%%%%%%%%%%%%%%%%%%%%%%%%%%%%%%%%%%%%%%%%%%%%%%%%%%%%%%%%%%%%%%%%%%%%%%%%%%%%%%%%%%%%%%%%%%%%%%%%%%%%%%%%%%%
\subsection{Graphs with edge weights} \label{3.1s:edges}
%%%%%%%%%%%%%%%%%%%%%%%%%%%%%%%%%%%%%%%%%%%%%%%%%%%%%%%%%%%%%%%%%%%%%%%%%%%%%%%%%%%%%%%%%%%%%%%%%%%%%%%%%%%%%
We first consider the case of a complete graph with $n$ nodes and real valued edge weights $g_{\{i,j\}}$. 
For example, the nodes of the graph could represent a set of points on a circle, and the weights between the 
nodes could be taken as the Euclidean distance between the corresponding points.

Denote by  ${\mathcal D}_g (G)$ the distribution of the weights of the graph.
Obviously, most graphs are not reconstructible from the distributions of their weights.
But in addition to the weights of the graph, 
one can also consider the sums of weights assigned to adjacent edges,  
which we denote by $\alpha_{i,j,k}$,
where
\[
\alpha_{i,j,k}=g_{\{i,j\}}+g_{\{j,k\}}.
\]
Denote by  ${\mathcal D}_\alpha (G)$ the distribution of the $\alpha_{i,j,k}$'s of a graph $G$. 
Observe that both  ${\mathcal D}_g (G)$ and  ${\mathcal D}_\alpha (G)$
are unchanged under a relabeling of the nodes.
Observe also that $\alpha_{i,j,k}=\alpha_{k,j,i}$, therefore
\begin{eqnarray*}
 \{ \alpha_{i,j,k} | i,j,k \text{ are distinct } \}&=& 
\{ \alpha_{i,j,k} | i,k\neq j, i<k \} \cup \{ \alpha_{i,j,k} | i,k\neq j, i>k \}, \\
&=& \{ \alpha_{i,j,k} | i,k\neq j, i<k \} \cup \{ \alpha_{k,j,i} | i,k\neq j, i>k \}, \\
&=&  \{ \alpha_{i,j,k} | i,k\neq j, i<k \} \cup  \{ \alpha_{i,j,k} | i,k\neq j, i<k \},
\end{eqnarray*}
and so to compute ${\mathcal D}_a (G)$, it is sufficient
to compute the distribution of the $\alpha_{i,j,k}$'s 
with $i,j,k$ distinct and $i<k$.
We now show that a large number of weighted graphs
are reconstructible from the distribution of their weights together 
with the distribution of the sum of weights assigned to adjacent edges. 

\begin{figure}%[t]%[htpb]
\begin{center}
Two Point Configurations on a Circle\\
\begin{tabular}{cc}
 \includegraphics[height=3.2cm]{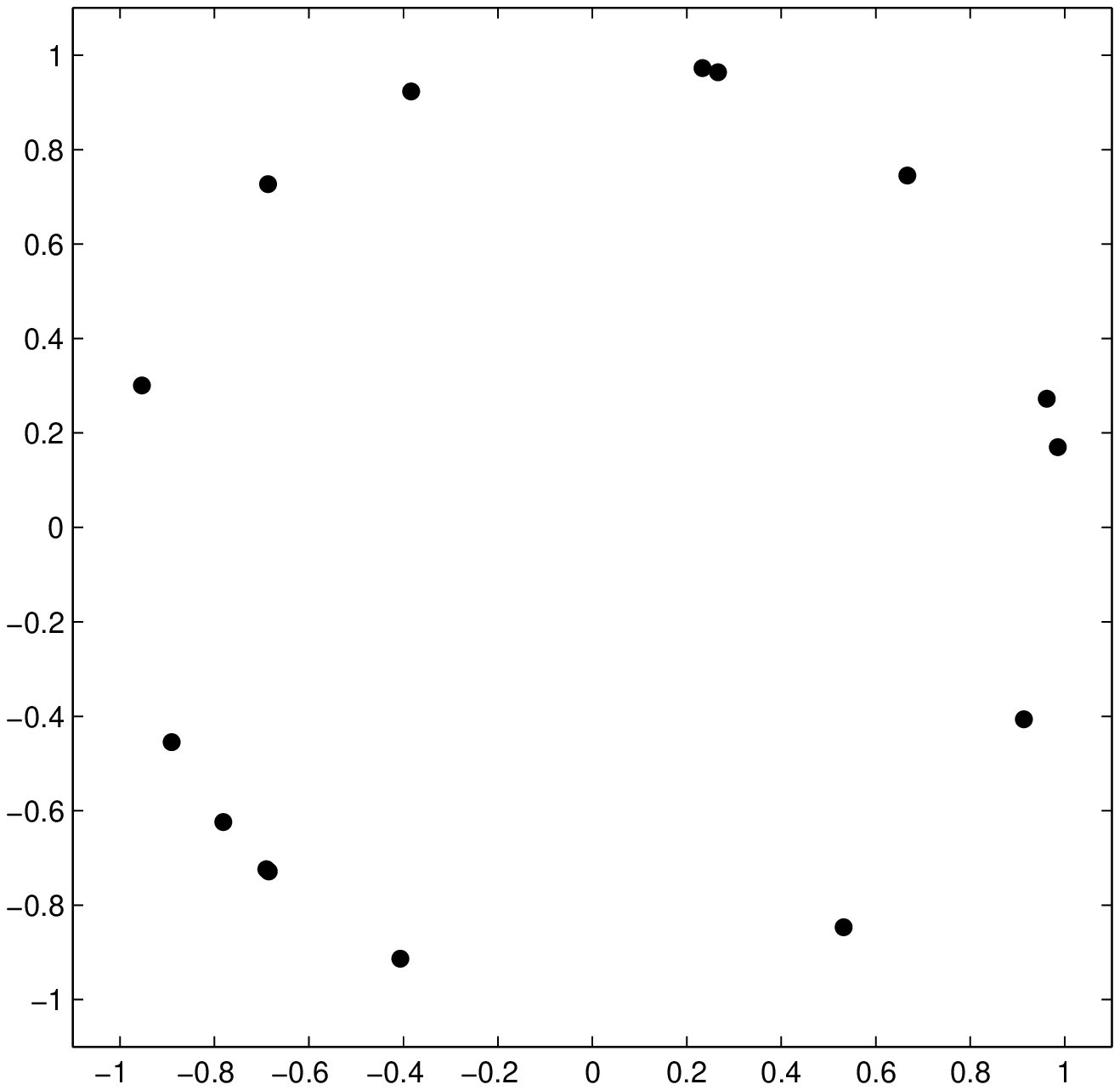} & \includegraphics[height=3.2cm]{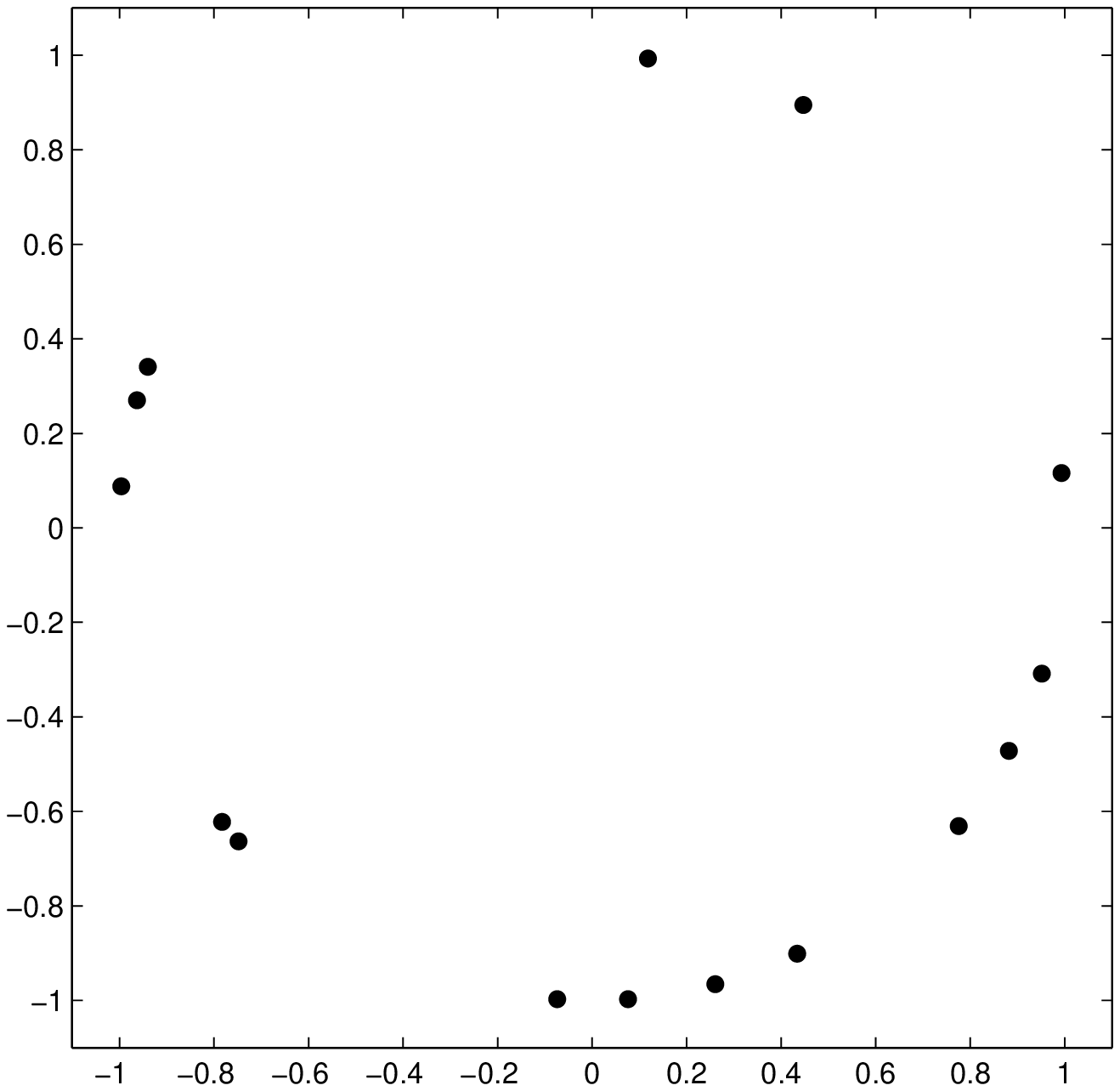} \\
\end{tabular}

Respective Histograms of Pairwise Distances \\
\begin{tabular}{cc}
 \includegraphics[height=3.2cm]{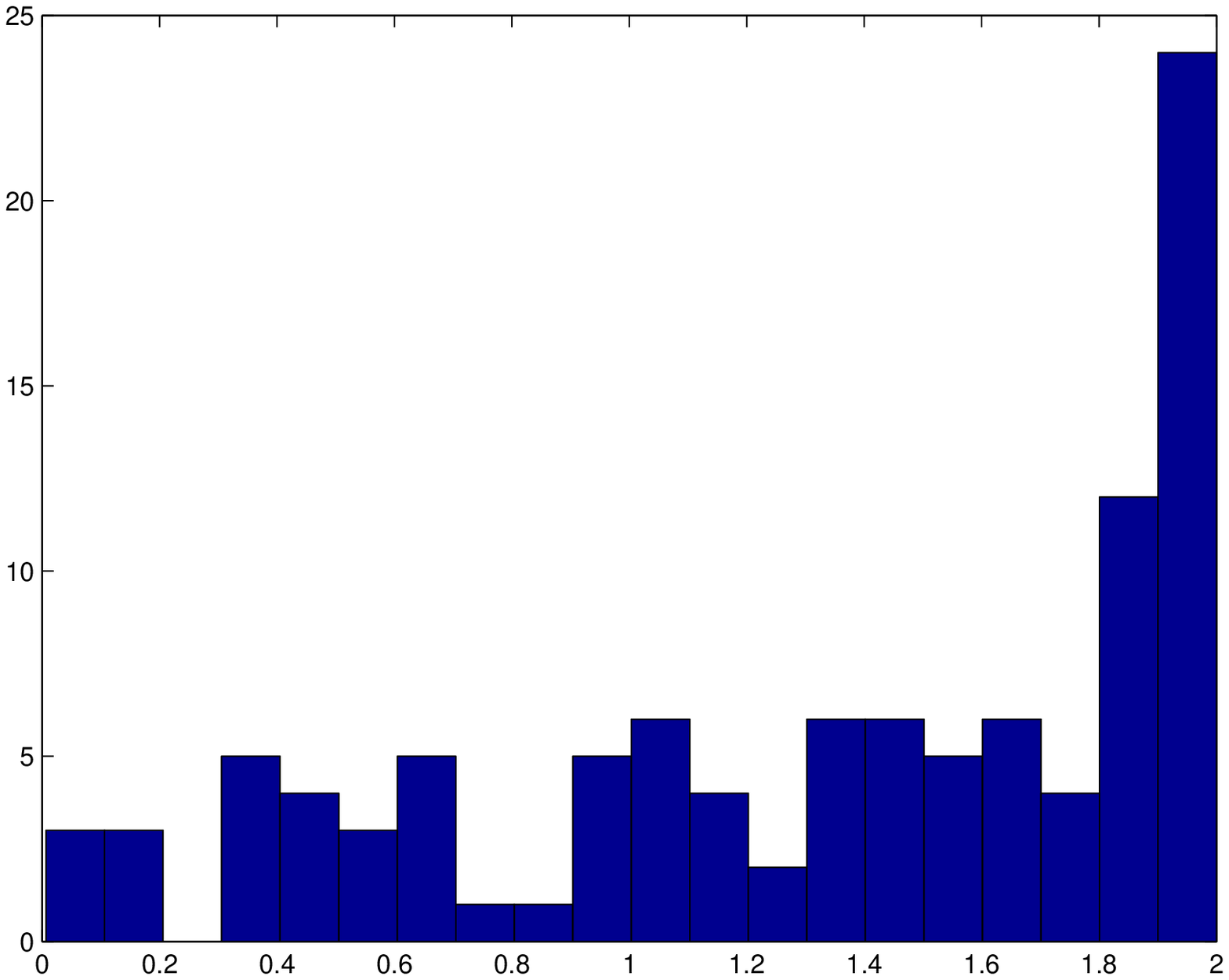} & \includegraphics[height=3.2cm]{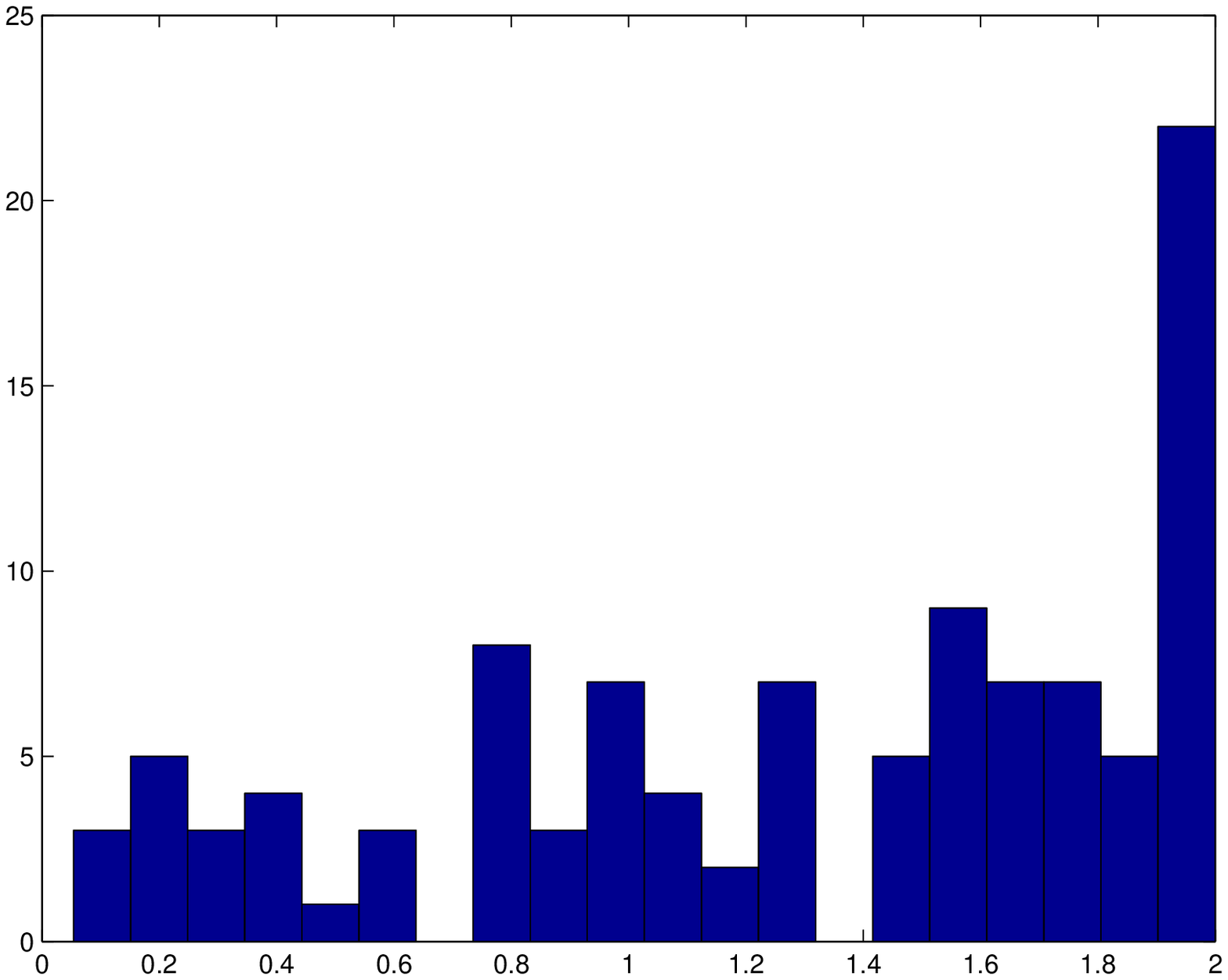} \\
\end{tabular}

Respective Histogram of Sums of Adjacent Distances  \\
\begin{tabular}{cc}
 \includegraphics[height=3.2cm]{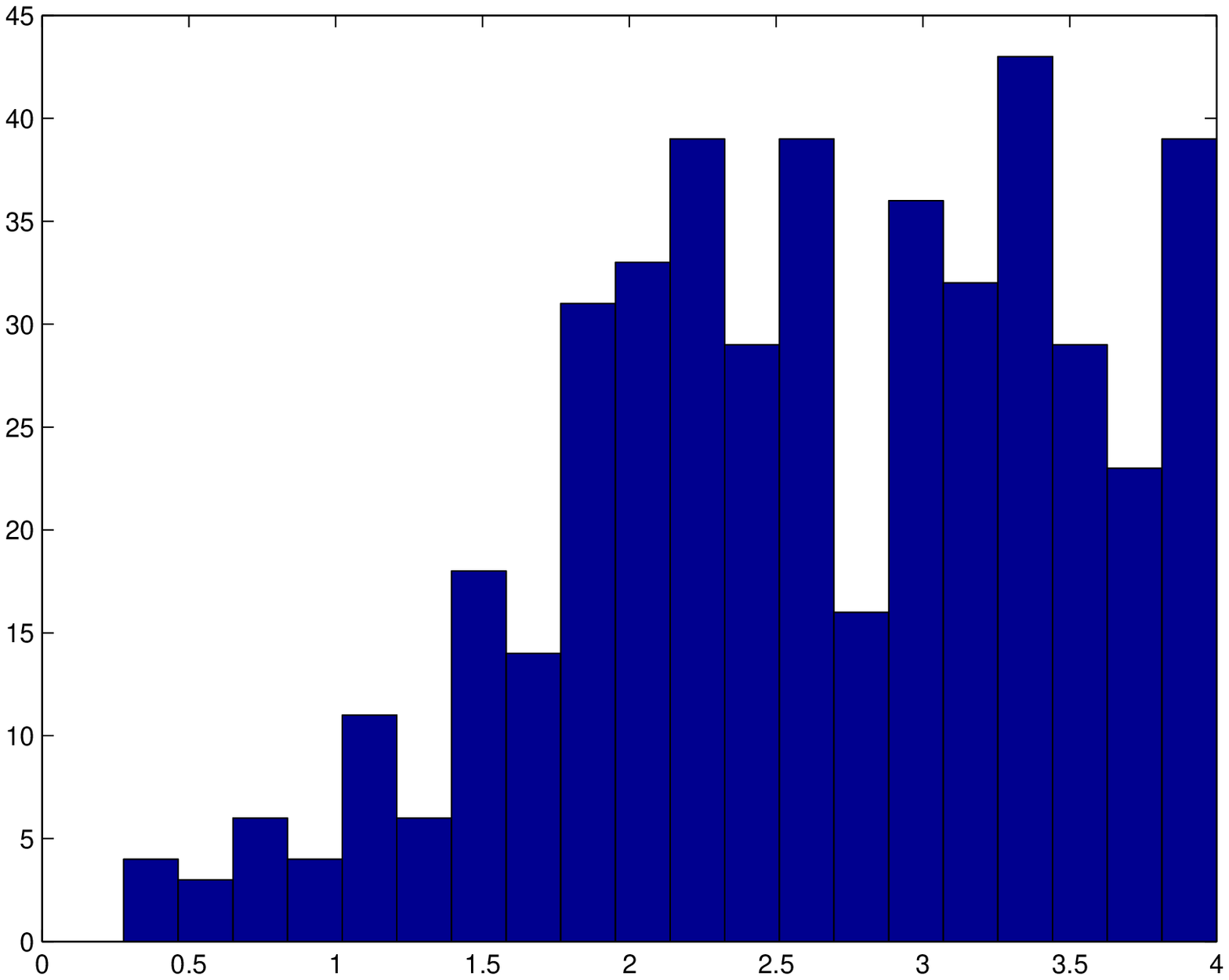} & \includegraphics[height=3.2cm]{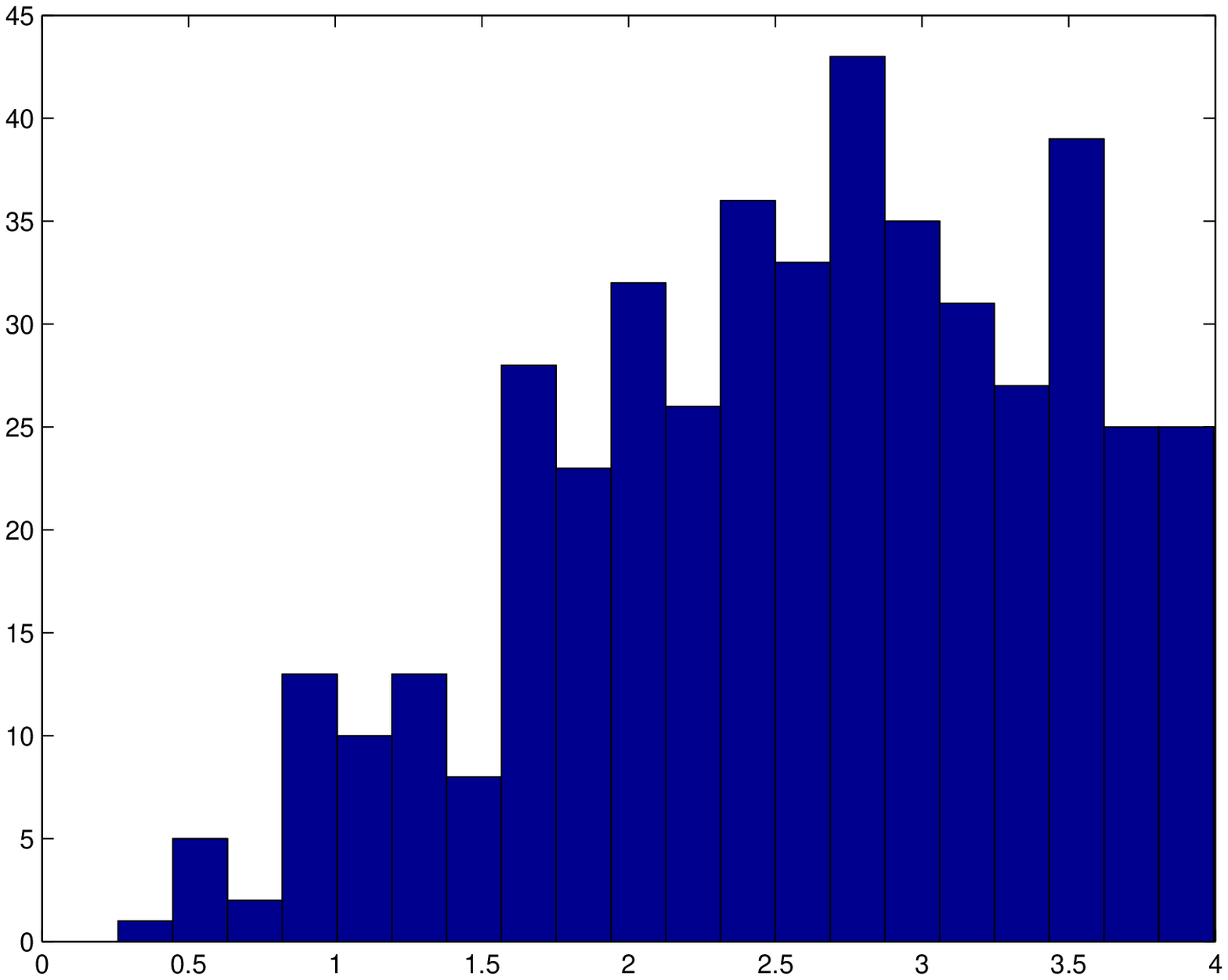} \\
\end{tabular}
\caption{A point configuration on a circle can be viewed as a complete weighted graph where the nodes correspond to the points and the edge weights are the Euclidean distances between the corresponding points. We propose to represent such a graph using two distributions: the distribution of the weights, and the distribution of the sums of any two adjacent weights. For most graphs, including the two point configurations pictured above, this is a lossless representation. More precisely, the set of graphs which are not uniquely reconstructible, up to isomorphism, from these two distributions for a set of measure zero. In particular, randomly chosen weights yield, with probability one, graphs that do not lie in the exceptional set. The points of the two configurations above were chosen (uniformly) at random on a unit circle; the histogram of their pairwise distances and the histogram of their sums of adjacent distances are clearly different.}
\end{center}
\label{circle}
\end{figure}

%%%%%%%%%%%%%%%%%%%%%%%%%%%%%%%%%%%%%%%%%%%
\begin{theorem}
\label{weighted graphs}
 Let  $G$ be a weighted graph with $n\geq 5$ nodes and weights $g_{\{i,j \}}\in {\mathbb R}$.
 Suppose that 
 \[ g_{\{i,j \}}+g_{\{j,k \}}\neq g_{\{m,p \}}+g_{\{q,r\}}, \]
 for every pairwise distinct $i,j,k$ 
 and every pairwise distinct $m,p,q,r$.
 Then $G$ is reconstructible 
 from  ${\mathcal D}_g (G)$ and ${\mathcal D}_\alpha (G)$.
\end{theorem}
%%%%%%%%%%%%%%%%%%%%%%%%%%%%%%%%%%%%%%%%%%%

\begin{proof}
 Let $\overline{G}=\{ \overline{g}_{\{i,j \}} \}$ 
 be another weighted graph with $n$ nodes such that
 ${\mathcal D}_g (\overline{G})={\mathcal D}_g (G)$ 
 and ${\mathcal D}_\alpha (\overline{G})={\mathcal D}_\alpha (G)$.
 Since  ${\mathcal D}_g (\overline{G})={\mathcal D}_g (G)$, there exists
 $\varphi\in S_{\binom{n}{2}}$ such that 
 \[ g_{\varphi \left(\{ i,j \} \right)}=\overline{g}_{\{i,j \}}, 
                  \text{ for every distinct }i,j=1,\ldots,n.\]
 We claim that there exists $\pi\in S_n$ such that 
 \[ \varphi \left( \{i,j \} \right)=\{\pi (i), \pi(j) \},
                   \text{ for every distinct }i,j=1,\ldots,n,\]
 that is to say, that $\varphi$ is simply a relabeling of the nodes.
 This is because, if we assume the contrary, 
 then, by \lref{L:permutations},
 there exists distinct indices $i_0,j_0,k_0$ such that
 \[\varphi\left(\{i_0,j_0\}\right)\cap
   \varphi\left(\{ j_0,k_0\}\right)=\emptyset.\]
 Since  ${\mathcal D}_\alpha (\overline{G})={\mathcal D}_\alpha (G)$, 
 there exists $i_1,j_1,k_1$ such that 
 \[ g_{\{ i_1,j_1\} }+g_{\{ j_1,k_1 \}}= 
         \overline{g}_{\{ i_0,j_0 \} }+\overline{g}_{\{ j_0,k_0\} }.\]
 But 
 \begin{eqnarray*}
  \overline{g}_{\{ i_0,j_0 \} }+\overline{g}_{\{ j_0,k_0\} } & =& 
                 {g}_{\varphi\left( \{ i_0,j_0\}\right) }
                       +{g}_{\varphi\left(\{ j_0,k_0\}\right)}\\
  &=& g_{\{m,p \}}+g_{\{ q,r\}},\text{ with }m,p,q,r\text{ distinct.}
 \end{eqnarray*}
 Therefore  
 $ g_{\{ i_1,j_1\} }+g_{\{ j_1,k_1 \}}= g_{\{m,p \}}+g_{\{ q,r\}}$ 
 which contradicts our hypothesis, 
 and thus $\varphi$ must be a relabeling of the nodes.
\end{proof}

\begin{remark}
A similar result would hold if we defined $\alpha_{i,j,k}$ as the product $\alpha_{i,j,k}=g_{\{i,j \}} g_{\{ j,k \}}$. 
In that case, the hypothesis would be that the graph must satisfy
\[ 
g_{\{i,j \}} g_{\{j,k \}}\neq g_{\{m,p \}} g_{\{q,r\}}, \]
 for every pairwise distinct $i,j,k$ 
 and every pairwise distinct $m,p,q,r$.
In particular, graphs with zero-valued weights would automatically be excluded. 
A work-around would be to shift all the weights by the same non-zero constant $\lambda$.
Alternatively, one could define  $\alpha_{i,j,k}$ 
as $\alpha_{i,j,k}=\left( g_{\{i,j \}} +\lambda\right) \left( g_{\{ j,k \}} +\lambda\right)$. 
Actually, the number of possibilities for the definition of $\alpha_{i,j,k}$ is endless, as any function $f$
of the two arguments $g_{\{i,j \}} g_{\{ j,k \}}$ can be used.
For any such function $f$, the corresponding theorem hypothesis would be written as
\[ 
f\left( g_{\{i,j \}} , g_{\{j,k \}}\right) \neq f\left( g_{\{m,p \}},  g_{\{q,r\}}\right), \]
 for every pairwise distinct $i,j,k$ 
 and every pairwise distinct $m,p,q,r$.
Similar remarks can be made for every theorem in this section.
\end{remark}

Observe that the proof of \tref{weighted graphs} is also valid
in the case of vector valued weights  ${\bf g}_{\{i,j \}}=\left( g^1_{\{i,j \}},\ldots, g^k_{\{i,j \}}  \right) \in \mathbb R^k$, for any integer $k$.
However, in that case, 
both ${\mathcal D}_g$ and ${\mathcal D}_\alpha$ become  $k$-dimensional distributions.
From a practical perspective, it is easier to deal with one-dimensional distributions.  
In particular, comparing one-dimensional distributions is much easier than comparing higher-dimensional distributions.
So an interesting question is: "Can graphs with edge weights ${\bf g}_{\{ i,j \}}\in {\mathbb R}^k$ 
be uniquely represented, up to isomorphism, by a set of one-dimensional distributions?".

We begin by considering the case ${\bf g}_{\{i,j \}}\in \mathbb
R^2$. (Such graphs will be called $(2,0)$-attribute graphs.)
Denote by  $\alpha^1_{i,j,k}$ the sum of the first weight assigned to two adjacent edges
and denote by  $\alpha^2_{i,j,k}$ the sum of the second weight assigned to same two adjacent edges:
\begin{eqnarray*}
        \alpha^1_{i,j,k} & = & g^1_{\{i,j \}}+g^1_{\{j,k \}},\\
        \alpha^2_{i,j,k} & = & g^2_{\{i,j \}}+g^2_{\{j,k \}},
\end{eqnarray*}
for all pairwise distinct $i,j,k$. We also define the mixed sum
\[
        \alpha^{12}_{i,j,k} =  g^1_{\{i,j \}}+g^2_{\{j,k \}}.
\]
Denote by  ${\mathcal D}_{\alpha^1} (G)$ and  ${\mathcal D}_{\alpha^2} (G)$  
the distributions of the
$\alpha_{i,j,k}^1$ and of the $\alpha_{i,j,k}^2$ respectively, 
for all $i,j,k$ pairwise distinct.
Similarly, denote by  ${\mathcal D}_{\alpha^{12}} (G)$ 
the distribution of the
$\alpha_{i,j,k}^{12}$ with $i,j,k$ pairwise distinct.
We show that {\em most } $(2,0)$-attribute graphs
are reconstructible from five distributions.

%%%%%%%%%%%%%%%%%%%%%%%%%%%%%%%%%%%%%%%%%%%
\begin{theorem}
 \label{2 edge attribute graphs}
 Let  $G$ be a graph with $n\geq 5$ nodes and vector valued weights ${\bf g}_{ \{ i,j \} }\in \RR^2$.
 Suppose that for every pairwise distinct $i,j,k$ 
 and every pairwise distinct $m,p,q,r$ we have
 \begin{eqnarray*}
 g^1_{\{i,j \}}+g^1_{\{j,k \}}\neq g^1_{\{m,p \}}+g^1_{\{q,r\}},\\
 g^2_{\{i,j \}}+g^2_{\{j,k \}}\neq g^2_{\{m,p \}}+g^2_{\{q,r\}},\\
 g^1_{\{i,j \}}+g^2_{\{j,k \}}\neq g^1_{\{m,p \}}+g^2_{\{q,r\}}.
 \end{eqnarray*}
 Then $G$ is reconstructible 
 from  ${\mathcal D}_{g^1} (G)$, 
 ${\mathcal D}_{g^2} (G)$, 
 ${\mathcal D}_{\alpha^1} (G)$, 
 ${\mathcal D}_{\alpha^2} (G)$ and ${\mathcal D}_{\alpha^{12}} (G)$.
\end{theorem}
%%%%%%%%%%%%%%%%%%%%%%%%%%%%%%%%%%%%%%%%%%%

\begin{proof}
 Let 
 $\overline{G}$ be another Graph with $n$ nodes and weights $\overline {\bf g}_{ \{ i,j \}} \in \RR^2$
 such that  ${\mathcal D}_{g^1} (\overline{G})={\mathcal D}_{g^1} (G)$, 
 ${\mathcal D}_{g^2} (\overline{G})={\mathcal D}_{g^2} (G)$, 
 ${\mathcal D}_{\alpha^1} (\overline{G})={\mathcal D}_{\alpha^1} (G)$, 
 ${\mathcal D}_{\alpha^2} (\overline{G})={\mathcal D}_{\alpha^2} (G)$ 
 and 
 ${\mathcal D}_{\alpha^{12}} (\overline{G})={\mathcal D}_{\alpha^{12}} (G)$.
 Since  ${\mathcal D}_{g^1} (\overline{G})={\mathcal D}_{g^1} (G)$,
 there exists $\varphi_1\in S_{\binom{n}{2}}$ such that
 \[ g^1_{\varphi_1 \left( \{ i, j\} \right) }=
                 \overline{g}^1_{\{ i,j\}}, \text{ for all distinct }
                      i,j=1,\ldots,n. \]
 By the same argument as in the proof of \tref{weighted graphs}, 
 there exists $\pi_1\in S_n$ such that 
 \[ \varphi_1 \left( \{ i, j\} \right)  = 
           \{ \pi_1(i), \pi_1(j)\},\text{ for all distinct }
                   i,j=1,\ldots,n.  \]
 Similarly, we can show that there exists $\pi_2\in S_n$ such that
 \[ g^2_{ \{ \pi_2(i), \pi_2(j)\} }=\overline{g}^2_{\{ i,j\}}, 
       \text{ for all distinct }
    i,j=1,\ldots,n. \]
 We claim that $\pi_1=\pi_2$.
 This is because, for any distinct $i,j,k\in \{1,\ldots,n \}$, 
 we have
 \begin{eqnarray*}
 \overline{\alpha}^{12}_{i,j,k}&=&\overline{g}^1_{\{ i,j \}}+\overline{g}^2_{\{ j,k\}}\\
 &=& g^1_{\{ \pi_1(i), \pi_1(j) \}}+g^2_{\{ \pi_2(j),\pi_2(k) \}}.
 \end{eqnarray*}
 But since ${\mathcal D}_{\alpha^{12}} (\overline{G})={\mathcal D}_{\alpha^{12}} (G)$,
 then for any $i,j,k \in \{ 1,\ldots,n\}$ distinct, 
 there exists $i',j',k'$ distinct
 such that 
 $\overline{\alpha}^{12}_{ijk}= \alpha^{12}_{i'j'k'}$.
 Therefore $ g^1_{\{ \pi_1(i), \pi_1(j) \} }+g^2_{\{ \pi_2(j),\pi_2(k) \}}=
 g^1_{\{ i',j' \}}+g^2_{\{ j',k' \}}$.
 By hypothesis, this implies that
 $\{\pi_1 (i), \pi_1 (j) \}\cap \{ \pi_2 (j), \pi_2 (k)\}\neq\emptyset$, 
 for every $i,j,k\in \{ 1,\ldots,n\}$ distinct.
 Let us choose three distinct indices $i_1,i_2,i_3$ 
 which are pairwise distinct from $j$ and $k$. 
 (We can do this because $n\geq 5$.)
 We have
 \begin{eqnarray*}
  (1)\phantom{bla}\{\pi_1 (i_1), \pi_1 (j) \}\cap \{ \pi_2 (j), \pi_2 (k)\} 
                               &\neq &\emptyset, \\
  (2)\phantom{bla} \{\pi_1 (i_2), \pi_1 (j) \}\cap \{ \pi_2 (j), \pi_2 (k)\} 
                             &\neq &\emptyset, \\
  (3)\phantom{bla} \{\pi_1 (i_3), \pi_1 (j) \}\cap \{ \pi_2 (j), \pi_2 (k)\} 
                              &\neq &\emptyset. 
 \end{eqnarray*}

 Assume that $\pi_1 (j)\notin\{ \pi_2 (j), \pi_2(k)\}$. 
 Then this means that, for all $l=1,2,3$, we have $\pi_1 (i_l)\in \{  \pi_2 (j), \pi_2 (k)  \}$,
 which contradicts the injectiveness of $\pi_1$. 
 We thus conclude that $\pi_1 (j)\in\{\pi_2(j),\pi_2(k) \}$, 
 for every distinct $j,k\in\{1,\ldots,n\}$.
 By varying the $k$, we obtain that $\pi_2(j)=\pi_1(j)$, for every $j$.
 \end{proof}

The above proof can trivially be generalized 
to the case of graphs with vector valued edge weights ${\bf g}_{ \{ i,j \} }\in \RR^d$, for any integer $d$.
Indeed, for any r $i,j,k$ pairwise distinct, we can define the quantities 
\begin{eqnarray*}
\alpha^l_{i,j,k} &=& g^l_{\{i,j \}}+g^l_{\{j,k \}},\text{ for all }l=1,\ldots,d, \\
\text{and }\alpha^{l,l+1}_{i,j,k} &=&g^l_{\{i,j \}}+g^{l+1}_{\{j,k \}},
\text{ for all }l=1,\ldots,d-1. 
\end{eqnarray*}
Then we denote by  ${\mathcal D}_{\alpha^l} (G)$ 
the distributions of the
$\alpha_{i,j,k}^l$, for all pairwise distinct $i,j,k$.
Similarly, denote by  ${\mathcal D}_{\alpha^{l,l+1}} (G)$ 
the distribution of the
$\alpha_{i,j,k}^{l,l+1}$, for all  pairwise distinct $i,j,k$.
Using the exact same arguments as for the above proof,
 we can show that these $3d-1$ one-dimensional distributions
fully characterize a large number of vector-valued weighted graphs. 
More precisely, we obtain the following theorem.

%%%%%%%%%%%%%%%%%%%%%%%%%%%%%%%%%%%%%%%%%%%
\begin{theorem}
 \label{d edge attribute graphs}
 Let  
 $G$ 
 be a graph 
 with $n\geq 5$ nodes and weights  ${\bf g}_{ \{ i,j \} }\in \RR^d$.
 Suppose that for every pairwise distinct $i,j,k$ 
 and every pairwise distinct $m,p,q,r$ we have
 \begin{eqnarray*}
 g^l_{\{i,j \}}+g^l_{\{j,k \}}&\neq& g^l_{\{m,p \}}+g^l_{\{q,r\}},
\text{ for all }l=1,\ldots,d,\\
\text{ and } g^l_{\{i,j \}}+g^{l+1}_{\{j,k \}}&\neq& g^l_{\{m,p \}}+g^{l+1}_{\{q,r\}},
\text{ for all }l=1,\ldots,d-1.
 \end{eqnarray*}
 Then $G$ is reconstructible 
 from the following $3d-1$ distributions:
\begin{eqnarray*}
 {\mathcal D}_{g^l} (G), & \text{ for all }&l=1,\ldots,d,\\
 {\mathcal D}_{\alpha^l} (G), &\text{ for all }&l=1,\ldots,d,\\
 {\mathcal D}_{\alpha^{l,l+1}} (G), &\text{ for all }& l=1,\ldots,d-1.
\end{eqnarray*}
\end{theorem}
%%%%%%%%%%%%%%%%%%%%%%%%%%%%%%%%%%%%%%%%%%%

%%%%%%%%%%%%%%%%%%%%%%%%%%%%%%%%%%%%%%%%%%%%%%%%%%%%%%%%%%%%%%%%%%%%%%%%%%%%%%%%%%%%%%%%%%%%%%%%%%%%%%%%%%%%%
\subsection{Graphs with node weights} \label{3.2s:nodes}
%%%%%%%%%%%%%%%%%%%%%%%%%%%%%%%%%%%%%%%%%%%%%%%%%%%%%%%%%%%%%%%%%%%%%%%%%%%%%%%%%%%%%%%%%%%%%%%%%%%%%%%%%%%%%

We now consider the case of a graph with $n$ nodes and node weights ${\bf w}_j={\bf g}_{\{ jj\}}\in \RR^d$.
First, let us assume that the graph does not have edge weights. 
Clearly, the graph is then reconstructible from the distribution of the node weights, 
which is a $d$-dimensional distribution.
But since we are interested in graph representation in terms of one-dimensional distributions,
we seek a different representation for $d>1$.

For any pairwise distinct $i,j,k$ and any $l\in\{1,\ldots,k-1\}$,
denote by $\beta^{l,l+1}_{i,j,k}$ the sum 
\[ \beta^{l,l+1}_{i,j,k}= w^l_i+ w^l_j+w^{l+1}_j + w^{l+1}_k.\]
Let ${\mathcal D}_{w^l} (G)$ be the distribution of the
$w_i^l$, for all $i=1,\ldots,n$.
Let  ${\mathcal D}_{\beta^{l,l+1}} (G)$ be the distribution of the
$b_{i,j,k}^{l,l+1}$, for all pairwise distinct $i,j,k$.

%%%%%%%%%%%%%%%%%%%%%%%%%%%%%%%%%%%%%%%%%%%
\begin{theorem}
 \label{k node attribute graphs}
 Let $G$  be a graph with $n\geq 5$ nodes and node weights ${\bf w}_j \in \RR^d$.
Suppose that, for every pairwise distinct $i,j,k$ 
 and every pairwise distinct $m,p,q,r$, we have
 \[  w^l_i w^l_j +w^{l+1}_j w^{l+1}_k \neq 
      w^l_m w^l_p+w^{l+1}_q w^{l+1}_r,\text{ for all }l=1,\ldots,d-1. \]
 Then $G$ is reconstructible from the following $2d-1$ distributions:
\begin{eqnarray*}
 {\mathcal D}_{w^l} (G),&\text{ for }& l=1,\ldots,d,\\
 {\mathcal D}_{\beta^{l,l+1}} (G), &\text{ for }& l=1,\ldots,d-1.
\end{eqnarray*}
\end{theorem}
%%%%%%%%%%%%%%%%%%%%%%%%%%%%%%%%%%%%%%%%%%%

\begin{proof}
 Let 
 $\overline{G}$ be another Graph with $n$ nodes and node weights  $\overline{{\bf w}}_j \in \RR^d$.
 such that  
\begin{eqnarray*}
{\mathcal D}_{w^l} (\overline{G}) &=& {\mathcal D}_{w^l} (G),\text{ for }l=1,\ldots, d,\\
\text{ and } 
{\mathcal D}_{\beta^{l,l+1}} (\overline{G}) &=& {\mathcal D}_{\beta^{l,l+1}} (G) ,\text{ for }l=1,\ldots, d-1. 
\end{eqnarray*}
 Since  ${\mathcal D}_{w^l} (\overline{G})={\mathcal D}_{w^l} (G)$,
 there exists $\pi_l \in S_n$ such that
 \[ w^l_{\pi (i)}=\overline{w}^l_i,\text{ for all }i=1,\ldots,n.\]
for $l=1,\ldots, d$.
. Using the same argument as in the proof of 
 Theorem\ref{2 edge attribute graphs} $d-1$ times,
 we can show that this implies that $\pi_1=\pi_2=\ldots \pi_d$.
Therefore $\overline{G}$ must be isomorphic to $G$.
\end{proof}

Finally, we consider the general case of a complete graph 
$G$ with both edge weights ${\bf g}_{\{ i,j\}}\in \RR^{d_1}$
and node weights  ${\bf w}_{i}\in \RR^{d_2}$.
For $l=1,\ldots,\min \{ d_1,d_2 \}$, 
let 
\[\Delta^l_{i,j,k}=w^l_i+ w^l_j+g^l_{\{j,k \}}.
\]
Denote by ${\mathcal D}_{\Delta^l}(G)$ the distribution
of the $\Delta^l_{i,j,k}$'s with $i,j,k$ pairwise distinct. 
Using the same arguments as for the previous theorems, 
we can show the following.

\begin{theorem}
Let $G$ be a graph with $n\geq 5$ nodes.
Suppose that $G$ has edge weights ${\bf g}_{\{ i,j\}}\in \RR^{d_1}$
and node weights  ${\bf w}_{i}\in \RR^{d_2}$.
Assume that 
\begin{eqnarray*}
%(1)
\phantom{bla} w^l_i+ w^l_j+w^{l+1}_j +w^{l+1}_k &\neq&  w^l_m + w^l_p+w^{l+1}_q + w^{l+1}_r, \text{ for all }l=1,\ldots,d_1, \\
%(2)
\phantom{bla} g^l_{\{ i,j\}} + g^{l+1}_{\{j,k\}}  
&\neq &  g^l_{\{ m,p\}} + g^{l+1}_{\{q,r\}}, \text{ for all }l=1,\ldots,d_2, \\
%(3)
\phantom{bla} w^1_i + w^1_j+g^1_{\{ j , k\}} &\neq &  
 w^1_m + w^1_p+g^1_{\{ q , r\}}, 
\end{eqnarray*}
for any pairwise distinct $i,j,k$ and any pairwise distinct $m,p,q,r$.
Then $G$ is reconstructible from the following distributions:
\begin{eqnarray*}
{\mathcal D}_{w^l} (G), &\text{ for }&l=1,\ldots,d_1, \\
{\mathcal D}_{\beta^{l,l+1}} (G), &\text{ for }&l=1,\ldots,d_1-1, \\
{\mathcal D}_{g^l} (G), &\text{ for }&l=1,\ldots,d_2, \\
{\mathcal D}_{\alpha^{l,l+1}} (G), &\text{ for }&l=1,\ldots,d_2-1, \\
\text{and }{\mathcal D}_{\Delta^1} (G).
\end{eqnarray*}
\end{theorem}

%%% Local Variables: 
%%% mode: latex
%%% TeX-master: "graphs"
%%% End: 

\bibliographystyle{mybibstyle} \bibliography{bib}

\bigskip

\begin{center}
\begin{tabular}{lll}
  Mireille Boutin & &Gregor Kemper \\
  School of Electrical and Computer Engineering  & & Technische Universit\"at M\"unchen \\
  Purdue University & & Zentrum Mathematik - M11 \\
  465 Northwestern Av. & & Boltzmannstr. 3 \\
  West Lafayette, IN & & 85\,748 Garching \\
  USA& & Germany \\
  {\tt mboutin$@$purdue.edu} & & {\tt kemper$@$ma.tum.de}
\end{tabular}
\end{center}

\end{document}